\newif\ifreport
\newif\ifarxiv

\arxivtrue
\reporttrue

\ifreport
	\ifarxiv
		\documentclass[referee,sn-mathphys-num]{sn-jnl}
	\else
		\documentclass[referee,lineno,sn-mathphys-num]{sn-jnl}
	\fi
\else
	\documentclass[iicol,sn-mathphys-num]{sn-jnl}
\fi

\usepackage[utf8]{inputenc} 
\usepackage{amssymb,amsmath}
\usepackage{amsfonts}
\usepackage{mathtools}
\usepackage{bbm}
\usepackage{cleveref}
\usepackage{graphicx}
\usepackage{subcaption}
\usepackage{csl}
\usepackage{comment}
\usepackage[normalem]{ulem} %required for strikeout
\usepackage[inline]{enumitem}

\newcommand{\titlestring}{Preserving Nonlinear Constraints in Variational Flow Filtering Data Assimilation}
\newcommand{\titlestringshort}{Preserving Nonlinear Constraints in Variational Flow Filters}
\newcommand{\authorstring}{Amit N. Subrahmanya, Andrey A. Popov, Reid J. Gomillion, Adrian Sandu}
\newcommand{\emailstring}{amitns@vt.edu, apopov@vt.edu, rjg18@vt.edu, sandu@cs.vt.edu}

\def\*#1{\boldsymbol{\mathbf{#1}}}
\def\!#1{\mathcal{#1}}

\newcommand{\syt}{\tau}

\newcommand{\x}{\*{x}}

\newcommand{\xf}{\*{x}_\mathrm{f}}
\newcommand{\xa}{\*{x}_\mathrm{a}}
\newcommand{\xt}{\*{x}_\syt}

\newcommand{\xfm}{\bar{\*{x}}_\mathrm{f}}

\newcommand{\xtm}{\bar{\*{x}}_\syt}

\newcommand{\xfk}{\*{x}_{\mathrm{f}, k}}
\newcommand{\xak}{\*{x}_{\mathrm{a}, k}}
\newcommand{\xakm}{\*{x}_{\mathrm{a}, k-1}}
\newcommand{\xtrk}{\*{x}_{\mathrm{true}, k}}

\newcommand{\xtro}{\*{x}_{\mathrm{true}, 0}}

\newcommand{\xai}{\*{x}_{\mathrm{a}, i}}
\newcommand{\xami}{\bar{\*{x}}_{\mathrm{a}, i}}
\newcommand{\xtri}{\*{x}_{\mathrm{true}, i}}

\newcommand{\X}{\*{X}}

\newcommand{\Xfk}{\*{X}_{\mathrm{f}, k}}
\newcommand{\Xak}{\*{X}_{\mathrm{a}, k}}
\newcommand{\Xakm}{\*{X}_{\mathrm{a}, k-1}}

\newcommand{\y}{\*{y}}
\newcommand{\z}{\*{z}}

\newcommand{\obserr}{\eta^\mathrm{o}}

\newcommand{\R}{\*{R}}

\newcommand{\Hn}{\!H}
\newcommand{\Mn}{\!M}
\newcommand{\Analy}{\!A}

\newcommand{\Hl}{\*H}

\newcommand{\dif}{\mathrm{d}}

\newcommand{\g}{\mathbf{g}}
\newcommand{\gx}{\mathbf{G}}

\newcommand{\Prob}{\mathcal{P}}
\newcommand{\Pa}{{\mathcal{P}_{\mathrm{a}}}}

\newcommand{\Pf}{{\mathcal{P}_{\mathrm{f}}}}
\newcommand{\Po}{{\mathcal{P}_{\mathrm{o}}}}
\newcommand{\Pt}{{\mathcal{P}_{\syt}}}

\newcommand{\CovPf}{\mathbf{P}_\mathrm{f}}

\newcommand{\CovPt}{\mathbf{P}_\syt}

\newcommand{\Rspace}{\mathbb{R}}

\newcommand{\nens}{n_\mathrm{e}}
\newcommand{\nstate}{n_\mathrm{s}}
\newcommand{\nobs}{n_\mathrm{o}}
\newcommand{\ncon}{n_\mathrm{c}}

\newtheorem{remark}{Remark}
\newtheorem{example}{Example}

\begin{document}

\ifreport
\ifarxiv
\csltitle{\titlestring}
\cslauthor{\authorstring}
\cslyear{24}
\cslreportnumber{1}
\cslemail{\emailstring}
\csltitlepage
\fi
\fi

\title[\titlestringshort]{\titlestring}

\author*[1]{\fnm{Amit N.} \sur{Subrahmanya}}\email{amitns@vt.edu}
\author[2]{\fnm{Andrey A.} \sur{Popov}}\email{apopov@vt.edu}
\author[1]{\fnm{Reid J.} \sur{Gomillion}}\email{rjg18@vt.edu}
\author[1]{\fnm{Adrian} \sur{Sandu}}\email{sandu@cs.vt.edu}

\affil*[1]{\orgdiv{Computational Science Laboratory, Department of Computer Science}, \orgname{Virginia Tech}, \orgaddress{\street{620 Drillfield Dr.}, \city{Blacksburg}, \postcode{24061}, \state{Virginia}, \country{USA}}}
\affil[2]{\orgdiv{Oden Institute for Computational Engineering and Sciences}, \orgname{University of Texas at Austin}, \orgaddress{\street{201 E. 24th Street}, \city{Austin}, \postcode{78712}, \state{Texas}, \country{USA}}}

\abstract{

Data assimilation aims to estimate the states of a dynamical system by optimally combining sparse and noisy observations of the physical system with uncertain forecasts produced by a computational model.
The states of many dynamical systems of interest obey nonlinear physical constraints, and the corresponding dynamics is confined to a certain sub-manifold of the state space.
Standard data assimilation techniques applied to such systems yield posterior states lying outside the manifold, violating the physical constraints. 

This work focuses on particle flow filters which use stochastic differential equations to evolve state samples from a prior distribution to samples from an observation-informed posterior distribution. 
The variational Fokker-Planck (VFP)---a generic particle flow filtering framework---is extended to incorporate non-linear, equality state constraints in the analysis. 
To this end, two algorithmic approaches that modify the VFP stochastic differential equation are discussed:
\begin{enumerate*}[label={(\roman*)}]
    \item VFPSTAB, to inexactly preserve constraints with the addition of a stabilizing drift term, and
    \item VFPDAE, to exactly preserve constraints by treating the VFP dynamics as a stochastic differential-algebraic equation (SDAE).
\end{enumerate*}
Additionally, an implicit-explicit time integrator is developed to evolve the VFPDAE dynamics.
The strength of the proposed approach for constraint preservation in data assimilation is demonstrated on three test problems: the double pendulum, Korteweg-de-Vries, and the incompressible Navier-Stokes equations.
}



\maketitle

%%%%%%%%%%%%%%%%%%%%%%%%%%%%%%%%%%%%%%%%%%%%%%%%%%%%%%%%%%%%%%%%%%%%%%%%
\section{Introduction}
\label{sec:intro}
%%%%%%%%%%%%%%%%%%%%%%%%%%%%%%%%%%%%%%%%%%%%%%%%%%%%%%%%%%%%%%%%%%%%%%%%

Data assimilation for state estimation attempts to combine uncertain, high-dimensional model simulations of a system with noisy (and typically low-dimensional) observations of the said system in a theoretically rigorous manner~\cite{Asch_2016_book,Reich_2015_book, Evensen_2022_book}. 
Data assimilation is widely used to improve forecasts in the areas of weather prediction, hydrology, seismology, and so on. 
Popular data assimilation methods include the well-known Kalman filter~\cite{Kalman_1960}, ensemble Kalman filter~\cite{Evensen_1994} (and its many variants such as the ensemble transform Kalman filter (ETKF)~\cite{Bishop_2001_ETKF}, ensemble adjustment Kalman filter (EAKF)~\cite{Anderson_2001_EAKF}, and the local ensemble transform Kalman filter(LETKF)~\cite{Hunt_2007_4DLETKF}), variational methods such as 3D-Var and 4D-Var ~\cite{Asch_2016_book,Evensen_2022_book} like, and particle filters~\cite{vanLeeuwen_2009_PF-review,vanLeeuwen_2019_PF-review}.

Particle flow filters solve the assimilation problem by iteratively evolving an ensemble of particles sampling the prior distribution into an ensemble of particles that sample the posterior through stochastic differential equations (SDE). 
While there are many flavors of particle flow filters~\cite{Pulido_2019_mapping-PF,Hu_2020_mapping-PF,Reich_2019_discrete-gradients,Reich_2021_FokkerPlanck,Stuart_2020_gradient-EnKF,Crouse_2019_consideration}, the main focus is on the authors'  \textit{Variational Fokker-Planck (VFP) approach ~\cite{Subrahmanya_2023_Ensemble}}, a flexible and generalized particle flow framework unifying many well-known particle flow methods. 
From the perspective of VFP, the particle flow in the state space is driven by a McKean-Vlasov-It\^{o} process whose drift minimizes the Kullback-Leibler divergence between the particle distribution and the posterior, and whose diffusion helps avoid particle collapse.
Multiple assumptions can be made on the prior and intermediate probability densities when computing the drift, and at the same time, can be combined with localization techniques to mitigate the curse of dimensionality~\cite{Hastie_2001_statsbook}.

Many physical systems of interest obey constraints such as conservation laws for mass, momentum, energy, enstrophy, etc~\cite{Arakawa_1972_Design,Janjic_2016_LETKF}.
In this work, the general term \textit{constraint} refers to any preserved fundamental property of a system. This work considers non-linear equality constraints. 
In general, while the computational models are built to preserve (in a numerical sense) the physical constraints of the system, standard data assimilation
algorithms are not. 
Traditional particle filters, which modify the weights associated with the particles (and not the particle states), preserve nonlinear constraints but are infeasible for high-dimensional problems~\cite{vanLeeuwen_2019_PF-review,Hastie_2001_statsbook}.
On the other hand, ensemble Kalman-like algorithms can conserve linear equality constraints~\cite{Janjic_2014_Preserve}, but not non-linear, or inequality constraints.
Like most standard filtering methods, the VFP framework~\cite{Subrahmanya_2023_Ensemble} makes no guarantees about constraint preservation.

Modifying filtering algorithms to achieve constraint preservation is not a new idea in data assimilation. Early work on preserving positivity constraints was done by Massicotte et. al~\cite{Massicotte_1995_Positivity}, where they \textit{inexactly} preserve positivity on filtered spectrometric data by heuristically damping the non-negativity of the state by multiplication with a deflation factor $\in [0, 1]$. 
Simon and Chia~\cite{Simon_2002_Equality} discuss Kalman filtering with state equality constraints.
They describe a projection-based formulation for linear equality constraints on the states and extend it to non-linear equality constraints by linearization.
They also present methods to restrict the Kalman gain such that the solution lies on the constraint manifold.
However, they conclude that the method of linearization may result in convergence issues when projecting onto the constraint manifold.

Julier and LaViola Jr.~\cite{Julier_2007_Nonlin} review and analyze two families of non-linear equality-constrained Kalman Filters, namely the pseudo-observation and projection methods. 
They identify two types of constraints, one that acts on the entire distribution (Type 1) and the other that acts only on the mean of the distribution (Type 2).
They conclude that the computationally simpler pseudo-observation method
(which includes constraints as observations with heuristically chosen error covariances)
preserves neither constraint; while the more computationally expensive projection method can preserve both types of constraints.
Gupta and Hauser~\cite{Gupta_2007_Constraint} introduce methods for inequality and equality linear constraint by the projection method or by modifying the Kalman gain. 
Sircoulomb et. al.~\cite{Sircoulomb_2008_Ineq} present a generic framework for dealing with non-linear inequality constraints by iterative projection.
Zanetti et. al.\cite{Zanetti_2009_norm} modify the Kalman gain operator for quadratic constraints, an approach useful when the norm constraint considers a small subset of the state space. 
Prakash et. al.~\cite{Prakash_2010_Ensemble}, propose a method to incorporate box constraints by modifying the ensemble Kalman filter to work on a truncated Gaussian distribution through a constrained optimization.
Bavdekar et. al.~\cite{Bavdekar_2013_Constrained} extend the previous idea to dual state and parameter estimation using the ensemble Kalman filter.
Their methodology---to separate the state and (constrained) parameter estimation---resulted in a lower error and a lower variance for the parameter estimate at increased cost (coming from additional model evolutions) when compared with a joint (state and parameter) EnKF approach.

Janjic et. al.~\cite{Janjic_2014_Preserve} deal with positivity preserving constraints in the ensemble Kalman filter by formulating the update as a sequence of quadratic programming problems.
Additionally, they show that the ensemble Kalman filters always conserve linear equality constraints assuming a preconsistent forecast. 
However, they also state that localization prevents this conservation.
Zeng and Janjic~\cite{Janjic_2016_LETKF} analyze the effect of localization radii, observed variables, and observation sparsity w.r.t the discretized grid on the conservation of mass, energy, enstrophy, divergence, and noise for the two-dimensional shallow water equations.
They observed that a small localization radius or observing only the heights degraded the conservation properties.
They also observed that the enstrophy violation, which should be minimal for stable non-linearities~\cite{Arakawa_1972_Design,Kacimi_2013_Arakawa}, is highly correlated to the localization radius and surprisingly, the ensemble inflation to prevent filter divergence.
Zeng, Janjic, et. al.~\cite{Janjic_2017_Preserve} extend the quadratic penalty idea~\cite{Janjic_2014_Preserve} to LETKF.

Li et. al.~\cite{Li_2019_Constrained} impose inequality constraints by projecting the unconstrained posterior Gaussian analysis particles onto the constrained region such that the constrained distribution is close in KL divergence to the unconstrained distribution.
Albers et. al ~\cite{Albers_2019_Constraint} discuss imposing linear equality and inequality constraints into the ensemble Kalman filter, by constraining the underlying EnKF quadratic optimization problem when there is a constraint violation.
They also extend it to the ensemble Kalman Inversion methodology~\cite{Stuart_2017_EnKF-inversion} to solve generic inverse problems.
% Stuart_2013_EnKF-inversion, check citation
Chada et. al.~\cite{Neil_2019_Box} incorporate box constraints into the Ensemble Kalman Inversion (EKI) by using a projected gradient method~\cite{Bertsekas_1982_PN}, into the gradient flow structure of the continuous limit EKI. 
The ensemble Kalman inversion idea relates somewhat to VFP as both methods involve a gradient flow to move samples from a prior distribution to a posterior distribution. 
Amor et. al~\cite{Amor_2018_Constrained} present a comprehensive overview of constraint-preserving techniques for ensemble Kalman filters and particle filters. 

The new contributions of this work are as follows.
\begin{enumerate}
    \item The incorporation of constraint preservation methodology into the generic VFP framework~\cite{Subrahmanya_2023_Ensemble} to obtain physically consistent (w.r.t the constraint) states.
    \item The first method---\textbf{VFPSTAB}---stabilizes the particle flow drift with an additional term to ensure the states are ``close enough'' to the constraint manifold.
    VFPSTAB is computationally efficient but inexact for constraint preservation. 
    \item The second method---\textbf{VFPDAE}---adds a constraint preservation term along with an algebraic constraint to the flow to ensure the states ``exactly'' respect the constraint manifold at all times during the filtering process. 
    Both these methods lie in stark contrast with historical methods where the assimilation and constraint preservation steps are treated independently, potentially leading to situations where the posterior and the constraint cannot be reconciled. 
    \item Finally, an implicit-explicit time stepping scheme is developed to evolve VFPDAE flow in pseudo-time efficiently.
\end{enumerate}

The remainder of the paper is organized as follows.
\Cref{sec:daconstraint} describes the data assimilation problem, focusing on filtering, and then presents the notation used in this work.
\Cref{sec:cons} examines the role of constraints in data assimilation along with a description of the traditional projection and pseudo-observation constraint preservation methods.
\Cref{sec:vfp} introduces the variational Fokker-Planck method, and develops two constraint preserving extensions, specifically, the VFPSTAB in \Cref{subsec:VFPSTAB} and the VFPDAE in \Cref{subsec:vfpdae}.
Next, \Cref{sec:expts} describes techniques for regularized covariance estimation in \Cref{subsubsec:covreg} and then compares VFPDAE and VFPSTAB to modified versions of the ETKF, and LETKF using the double pendulum (\Cref{subsec:exptdp}), Korteweg-de Vries (\Cref{subsec:exptkdv}) and incompressible Navier-Stokes (\Cref{subsec:exptins}) problems. 
Finally, the conclusions are presented in \cref{sec:conc}.
%

%%%%%%%%%%%%%%%%%%%%%%%%%%%%%%%%%%%%%%%%%%%%%%%%%%%%%%%%%%%%%%%%%%%%%%%%
\section{Data assimilation}
\label{sec:daconstraint}
%%%%%%%%%%%%%%%%%%%%%%%%%%%%%%%%%%%%%%%%%%%%%%%%%%%%%%%%%%%%%%%%%%%%%%%%

A discrete-time dynamical system model is described by an operator that propagates the states:
\begin{equation}
\label{eqn:model-dynamics}
    \x_{k} = \Mn_{k}(\x_{k-1}), \quad k \ge 1,
\end{equation}
where $k$ indexes the time $t_k$, $\x_{k - 1} \in \Rspace^{\nstate}$ represents the states of the system at time $t_{k-1}$, and the model operator $\Mn_{k}: \Rspace^{\nstate} \to \Rspace^{\nstate}$ evolves the states from time $t_{k-1}$ to time $t_k$. 
This paper considers deterministic models $\Mn$ described by (partial) differential-algebraic equations~\cite{Ascher_1998_DAE}.

The model states must be combined with observations of the unknown reality to track the true states. 
These observations of an unknown true state are typically sparse (in space and time) and noisy, and given by 
\begin{equation}
    \y_{k} = \Hn_{k}(\xtrk) + \obserr_k,
\end{equation}
where $\y_k \in \mathbb{R}^{\nobs}$ is the observation at time $t_k$, $\Hn_{k}: \Rspace^{\nstate} \to \Rspace^{\nobs}$ the observation operator and $\obserr_k$ is the observation noise. 
In a typical data assimilation setting, $\nobs \ll \nstate$ as it is infeasible to have a larger set of observations. 

Filtering alternates between forecasting, where previous analysis states are propagated to a future time using the computational model; and analysis, which combines this forecast with the current observation to obtain new analysis states.
This cycle is summarized as
\begin{equation}\label{eq:basicda}
\begin{aligned}
    \xfk &= \Mn_{k}(\xakm), &\quad &\mathrm{(Forecast)} \\
    \xak &= \Analy(\xfk, \y_k; \*{\theta}), &\quad &\mathrm{(Analysis)}
\end{aligned} 
\end{equation}
where $\*{\theta} \in \mathbb{R}^p$ is a p-dimensional set of parameters of the analysis procedure.
For example, in the localized ensemble Kalman Filter, the parameters $\*{\theta}$ could include ensemble inflation and localization radius. 
Through the Bayes' rule, we have 
\begin{equation}
    \Prob(\xak) = \Prob(\xfk\mid\y_k) = \frac{\Prob(\y_k\mid\xfk)\cdot  \Prob(\xfk)}{\Prob(\y_k)},
\end{equation}
where the analysis distribution is the (Bayesian) posterior that combines the observational likelihood $\Prob(\y_k\mid\xfk)$ and the prior distribution of the forecast $\Prob(\xfk)$.
Equivalently, we refer to the probabilities as $\Prob(\xfk) \equiv \Pf(\x_k)$, $\Prob(\y_k \mid \xfk) \equiv \Po(\y_k - \Hn_{k}(\x_k))$, and $\Prob(\xak) \equiv \Pa(\x_k)$.

This paper exclusively considers ensemble approaches to filtering, i.e. 
\begin{equation}\label{eq:basiceda}
\begin{aligned}
    \Xfk &= \Mn_{k}(\Xakm), &\quad &\mathrm{(Forecast)} \\
    \Xak &= \Analy(\Xfk, \y_k; \*{\theta}), &\quad &\mathrm{(Analysis)}
\end{aligned} 
\end{equation}
where the ensemble
\begin{equation}
    \X_{\ast, k} = 
    \begin{bmatrix}
        \x_{\ast, k}^{[1]} & \x_{\ast, k}^{[2]} & \dots & \x_{\ast, k}^{[\nens]}  
    \end{bmatrix} \in \Rspace^{\nstate \times \nens},
\end{equation}
has $\nens$ exchangeable random variables, representing samples (not necessarily independent) from the same distribution.
Here, it is assumed that $\xfk^{[e]} \sim \Pf(\x_k)$ and $\xak^{[e]} \sim \Pa(\x_k)$ for all $e$ in the index set $\!I_{\nens} = \{ 1 , 2 , \dots , \nens \}$.
%

%%%%%%%%%%%%%%%%%%%%%%%%%%%%%%%%%%%%%%%%%%%%%%%%%%%%%%%%%%%%%%%%%%%%%%%%
\section{Constraints}
\label{sec:cons}
%%%%%%%%%%%%%%%%%%%%%%%%%%%%%%%%%%%%%%%%%%%%%%%%%%%%%%%%%%%%%%%%%%%%%%%%

Many dynamical systems of interest, possess natural constraints on their states~\cite{Simon_2002_Equality, Gupta_2007_Constraint, Janjic_2014_Preserve} as in \cref{eqn:model-invariants}.
For example, 
\begin{enumerate*}[label={(\roman*)}]
    \item geophysical systems obey conservation laws such as mass, momentum, energy, and so on;
    \item mechanical systems are constrained by lengths or velocities of its components.
\end{enumerate*}
These constraints, which describe physical or geometric laws, are ignored while filtering as they can be challenging to incorporate into the filtering methodology~\cite{Simon_2002_Equality}.

For spatially distributed dynamical systems (typically described by PDEs), constraints are classified as global and local~\cite{Anco_2020_Conservation}. 
Global constraints preserve information across the entire domain of the problem. 
For example, the conservation of total energy of a closed system, the conservation of enstrophy for a closed fluid simulation, the length of a pendulum, etc. 
Local constraints preserve information in a pointwise fashion (i.e. at spatial locations) across the domain. 
For example, this can be the continuity equation for an incompressible fluid where the divergence is pointwise zero at every spatial location across the domain; or that of non-negativity where the concentration of a chemical is greater than or equal to zero at each point across the domain. 
Constraints can also be divided into two types based on their origin.
\begin{enumerate}
    \item Constraints defined by state invariants that do not change with either time or across samples. 
    For models such as \cref{eqn:model-dynamics}, the constraints are described by the algebraic relations:
    \begin{equation}
    \label{eqn:model-invariants}
    \g(\x_{k-1})=\*0 \ \Rightarrow \ \g(\x_{k})=\g\left(\Mn_{k}(\x_{k-1})\right)=\*0,
    \end{equation}
    where $\g : \Rspace^{\nstate} \to \Rspace^{\ncon}$ is the constraint function.
    If $\x_{k - 1}$ were to be modified as $\tilde{\x}_{k - 1}$ such that $\g(\tilde{\x}_{k - 1}) \neq \*0$, the following situations could occur with $\tilde{\x}_k = \Mn_{k}(\tilde{\x}_{k - 1})$.
    \begin{enumerate}
        \item Some model dynamics naturally pull back the states onto the constraint manifold, i.e. $\g\left(\tilde{\x}_k\right) = \*0$.
        The states need no special treatment in this scenario and such models are not considered in this work.
        \item Some models permit a small degree of flexibility, i.e. 
        \begin{equation}
            |\g(\tilde{\x}_{k - 1})| \le \varepsilon \quad \Rightarrow \quad |\g(\tilde{\x}_k)| \le \varepsilon,
        \end{equation}
        where the constraints are \textit{weakly} preserved.
        \item Sometimes $\tilde{\x}_k$ is not computable if $\g(\tilde{\x}_{k - 1}) \neq \*0$ or even if $\tilde{\x}_k$ is computable, $\g(\tilde{\x}_k) \neq \*0$. 
        For these situations, $\tilde{\x}_{k - 1}$ must be modified to lie on the constraint manifold and then evolved.
    \end{enumerate}
    The following equations are all at time $t_k$ (subscript $k$ has been dropped in the equations for simplicity), the data assimilation problem assumes a preconsistent forecast i.e. 
    \begin{equation}
        \forall e \in \!I_{\nens} \quad \g(\xf^{[e]}) = \*0, \quad \text{or} \quad |\g(\xf^{[e]})| \le \varepsilon.
    \end{equation}
    The goal is to have analysis ensemble members who also live on the constraint manifold as
    \begin{equation}
        \forall e \in \!I_{\nens} \quad \g(\xa^{[e]}) = \*0, \quad \text{or} \quad |\g(\xa^{[e]})| \le \varepsilon.
    \end{equation}
    \begin{remark}
        The constraints can also vary across time, or differ across the ensemble.
        We do not delve into these cases, because the underlying ideas regarding these constraints and their preservation are mostly similar to the above.
    \end{remark}
    \item Constraints that preserve information between the forecast and analysis samples.
    Essentially, the constraints (again, all at time $t_k$) are defined as
    \begin{equation}
        \forall e \in \!I_{\nens} \quad \g^{[e]}(\xa^{[e]}; \xf^{[e]}) = \*h(\xa^{[e]}) - \*h(\xf^{[e]}),
    \end{equation}
    where $\*h :  \Rspace^{\nstate} \to \Rspace^{\ncon}$ represents some preservable quantities.
    The constraint preservation requires
    \begin{equation}
        \forall e \in \!I_{\nens} \quad \g^{[e]}(\xa^{[e]}; \xf^{[e]}) = \*0, \quad \text{or} \quad |\g^{[e]}(\xa^{[e]}; \xf^{[e]})| \le \varepsilon.
    \end{equation}
    Typically, these constraints are artificial (rather than natural) and enforced based on need. 
\end{enumerate}
Problems with both constraint types can be formulated, such as the incompressible Navier Stokes experiment in~\cref{subsec:exptins}.
With a few exceptions, most filtering methods ignore state constraints, and as a consequence, fail to preserve them~\cite{Simon_2002_Equality,Janjic_2016_LETKF}.

\begin{example}[Pendulum]
Consider an example problem with the simple pendulum to motivate constraint preservation further.
Consider the situation where the forecast pendulum ensemble, the truth, the observation, and the constraint manifold are as depicted in \Cref{fig:example_p}.
The constraint manifold here is defined by the geometric constraints of the pendulum such as its length and the extreme position (dictating the total mechanical energy of the system in the absence of forcing or dissipation).
While the ETKF analysis (labeled \textit{Standard ETKF} in \Cref{fig:example_p}) ensemble is closer to the observations, it is completely off the constraint manifold.
These inconsistent analysis states are meaningless as the pendulum's length must remain constant after filtering, and hence, a better and physically consistent analysis is preferred (labeled \textit{Constrained ETKF} in \Cref{fig:example_p}).
\begin{figure*}[!ht]
    \centering      
    \includegraphics[width=0.5\linewidth]{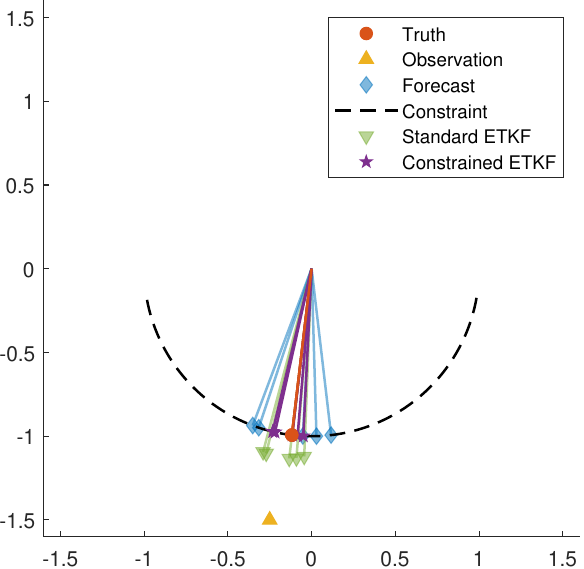}
    \caption{
    A simple pendulum example where the forecast ensemble members lie on the constraint manifold, but the standard ETKF analysis does not. The goal is to obtain an analysis ensemble that lives on the constraint manifold like the Constrained ETKF.} 
    \label{fig:example_p}
\end{figure*}
\end{example}

We now discuss two different approaches to constrain analyses.

%%%%%%%%%%%%%%%%%%%%%%%%%%%%%%%%%%%%%%%%%%%%%%%%%%%%%%%%%%%%%%%%%%%%%%%%
\subsection{Projection onto constraint manifold after filtering}
\label{subsec:proj}
%%%%%%%%%%%%%%%%%%%%%%%%%%%%%%%%%%%%%%%%%%%%%%%%%%%%%%%%%%%%%%%%%%%%%%%%

The projection after filtering approach has two main steps.
First, a constraint agnostic filtering scheme is applied to the forecast to obtain the unconstrained analysis: 
\begin{equation}
\widehat{\x}_{\rm a}^{[e]} = \Analy(\xf^{[e]}, \y, \*{\theta}),  \quad  \forall e \in \!I_{\nens}.    
\end{equation}
Next, each analysis ensemble member is projected onto the manifold using:
\begin{subequations}
\label{eqn:analysis-projection}
\begin{equation}
\label{eqn:analysis-projection-x}
    \xa^{[e]} = \widehat{\x}_{\rm a}^{[e]} - \gx^\top(\widehat{\x}_{\rm a}^{[e]})\, \*\lambda^{[e]} ,  \quad  \forall e \in \!I_{\nens},
\end{equation}
where $\*\lambda^{[e]} \in \Rspace^{\ncon}$ is a solution of the nonlinear equation:
\begin{equation}
\label{eqn:analysis-projection-lambda}
    \*0 =  \g\left(\xa^{[e]}\right) = \g\left(\widehat{\x}_{\rm a}^{[e]} - \gx^\top(\widehat{\x}_{\rm a}^{[e]})\, \*\lambda^{[e]} \right),  \quad  \forall e \in \!I_{\nens},
\end{equation}
\end{subequations}
where the Jacobian of the constraints defines the tangent constraint manifold:
\begin{equation}
\label{eqn:constraint-Jacobian}
\gx(\x) \coloneqq \frac{d\g(\x)}{d \x}.
\end{equation}
The approach \cref{eqn:analysis-projection} is to move the unconstrained analysis states along the adjoint constraint manifold $\gx^\top$ until it lands onto the constraint manifold.
This idea is similar to projected Runge-Kutta methods for differential-algebraic equations~\cite{Ascher_1998_DAE}. 
For linear constraints, formulation \cref{eqn:analysis-projection} leads to the Mean Square Method derivation of projection from Simon and Chia~\cite{Simon_2002_Equality}.

%%%%%%%%%%%%%%%%%%%%%%%%%%%%%%%%%%%%%%%%%%%%%%%%%%%%%%%%%%%%%%%%%%%%%%%%
\subsection{Pseudo-observations of constraints}
\label{subsec:pseudo}
%%%%%%%%%%%%%%%%%%%%%%%%%%%%%%%%%%%%%%%%%%%%%%%%%%%%%%%%%%%%%%%%%%%%%%%%

The pseudo-observation formulation here follows Julier and LaViola Jr~\cite{Julier_2007_Nonlin}.
Formally, the standard observation operator $\Hn(\x)$ is augmented with $\g(\x)$ to obtain the extended observation operator $\widehat{\Hn}$:
\begin{equation}
    \widehat{\Hn}(\x) \coloneqq \begin{bmatrix}
        \Hn(\x) \\ \g(\x)
    \end{bmatrix},
\end{equation}
and the data with a pseudo-observation of the constraint value, to obtain the extended observation $\widehat{\y}$:
\begin{equation}
    \widehat{\y} \coloneqq \begin{bmatrix}
        \y \\ \*0
    \end{bmatrix}.
\end{equation}
For any ensemble method, the observation error covariance $\R$ is augmented with a user-defined covariance $\R^{\g}$ that defines an acceptable amount of constraint violation: 
\begin{equation}
    \widehat{\R} \coloneqq 
    \begin{bmatrix}
        \R & \*0\\\*0 & \R^{\g}
    \end{bmatrix}.
\end{equation}
The extended observation triplet $\widehat{\Hn}(\x)$, $\widehat{\y}$, and $\widehat{\R}$ is used in the standard way in filtering methods.
However, combining pseudo-observations with localization is difficult and is not explored here.

%%%%%%%%%%%%%%%%%%%%%%%%%%%%%%%%%%%%%%%%%%%%%%%%%%%%%%%%%%%%%%%%%%%%%%%%
\section{Ensemble Variational Fokker Planck Filters}
\label{sec:vfp}
%%%%%%%%%%%%%%%%%%%%%%%%%%%%%%%%%%%%%%%%%%%%%%%%%%%%%%%%%%%%%%%%%%%%%%%%
%
The ensemble variational Fokker Planck method~\cite{Subrahmanya_2023_Ensemble} for data assimilation evolves a set of samples from the prior distribution to the posterior distribution along a defined stochastic differential equation~\cite{Subrahmanya_2023_Ensemble} in pseudo-time.
Specifically, the samples move under a McKean-Vlasov-It\^{o} process whose drift minimizes the Kullback-Leibler divergence between the probability density of the samples and the target posterior and whose diffusion is tuned to prevent filter divergence.
This framework treats the density as a "parameter" estimated from the samples which can be chosen arbitrarily as long as the log-gradient of the density can be computed exactly or estimated.
Formally, the McKean-Vlasov-It\^{o} process defined as 
\begin{equation}\label{eq:vfpito}
    \dif \xt = \*F(\syt, \xt, \Pt) \dif \syt + \boldsymbol{\sigma}(\syt, \xt, \Pt) \dif \*W,
\end{equation}
whose corresponding Fokker-Planck equation is 
\begin{equation}\label{eq:vfpfokkerplanck}
    \frac{\partial \Pt(\x)}{\partial \syt} = -\nabla_{\x} \cdot \left( \Pt(\x) \*F  - \nabla_{\x} \cdot \left(\frac{\Pt(\x) \boldsymbol{\sigma}\boldsymbol{\sigma}^\top}{2} \right) \right),
\end{equation}
where $\syt$ is the pseudo-time, $\xt$ is a sample, $\Pt$ is the density of the sample, $\*F$ is the KL divergence minimizing optimal drift, $\boldsymbol{\sigma}$ is the diffusion tensor and $\*W$ is the Wiener process and $(\nabla_{\x} \cdot ): \Rspace^{\nstate \times \nstate} \to \Rspace^{\nstate}$ is the divergence term.
The initial conditions are given as
\begin{equation}
    \syt = 0, \;\; \xt = \xf \; \textrm{ and } \; \Pt = \Pf(\x).
\end{equation}
When \cref{eq:vfpito} is evolved in pseudo-time, we obtain the analysis samples as
\begin{equation}
    \syt \to \infty, \;\; \xt = \xa \; \textrm{ and } \; \Pt = \Pa(\x).
\end{equation}
In practice, the integration is terminated at some finite pseudo-time based on a criterion evaluating the convergence of the samples to the posterior.
The optimal drift that minimizes the KL-divergence between the ensemble estimate and the posterior (derived in \cite{Subrahmanya_2023_Ensemble}) is given as 
\begin{equation}\label{eq:vfpoptimaldrift}
\begin{split}
    \*F(\syt, \x, \Pt) &= \*{A}(\x, \Pt) \left( \nabla_{\x} \log{\Pa(\x)} - \nabla_{\x} \log{\Pt(\x)} \right)  \\
    &+ \left(\frac{\boldsymbol{\sigma}\boldsymbol{\sigma}^\top}{2}\right) \nabla_{\x} \log{\Pt(\x)} + \nabla_{\x} \cdot \left(\frac{\boldsymbol{\sigma}\boldsymbol{\sigma}^\top}{2}\right), 
\end{split}
\end{equation}
where $\*{A}(\x, \Pt)$ can be any symmetric positive definite matrix. 
The term $\nabla_{\x} \log{\Pa(\x)}$ pushes the samples toward the posterior, $\nabla_{\x} \log{\Pt(\x)}$ attempts to keep the samples spread apart (deterministically), and the terms containing $\boldsymbol{\sigma}\boldsymbol{\sigma}^\top/2$ account for the stochastic diffusion.
For a more rigorous description of this framework, we refer the readers to our previous work~\cite{Subrahmanya_2023_Ensemble}.

The following simplifying assumptions are made for the remainder of the paper: 
\begin{enumerate}[label={(\roman*)}]
    \item Firstly, $\*{A}(\x, \Pt) = \*I_{\nstate}$.
    This is a trivial choice (theoretically and computationally) and further justification is discussed in~\cite{Subrahmanya_2023_Ensemble}.
    \item  
    Deriving the exact flow based on $\Pt$, $\Po$, and $\Pf$ would require the estimation of probability densities on manifolds.
    This idea is not explored in this work as the main focus is on constraint preservation in the VFP framework.
    Thus, for simplicity, the computational flow assumes the densities $\Pt$, $\Po$, and $\Pf$ to all be Gaussian.
    While the Gaussian assumption works in practice, VFP is not restricted by this and can choose any quantifiable gradient log density to derive the flow (discussed in \cite{Subrahmanya_2023_Ensemble}).
    It is more likely than not, that a Gaussian in $\Rspace^{\nstate}$ will not be Gaussian on the constraint manifold.
    But practically, the experiments show that this is a good enough approximation for the flow.
    \item     
    The diffusion is assumed to be some scaled square root (or Cholesky factor) of the forecast covariance estimate $\boldsymbol{\sigma}(\syt, \xt, \Pt) = \alpha \sqrt{\CovPf}$ with $\alpha$ being a problem specific tunable parameter.
    As $\boldsymbol{\sigma}$ is independent of $\xt$, we have $\left(\frac{\boldsymbol{\sigma}\boldsymbol{\sigma}^\top}{2}\right) = \alpha^2 \CovPf$ and $\nabla_{\x} \cdot \left(\frac{\boldsymbol{\sigma}\boldsymbol{\sigma}^\top}{2}\right) = \*0$
    Note that for $\*{\sigma} \in \mathbb{R}^{\nstate \times \ell}$, where $\ell < \nstate$, the diffusion dynamics are restricted to the forecast ensemble subspace~\cite{Inigo_2020_Langevin}.
\end{enumerate}

\begin{example}[Gaussian case]
Under the Gaussian assumptions, 
\begin{equation}\label{eq:gaussiandrift}
    \begin{split}
    \nabla_{\x} \log{\Pt(\x)} &= -\CovPt^{-1}\left( \x - \xtm \right), \;\mathrm{and} \\
        \nabla_{\x} \log{\Pa(\x)} &= -\CovPf^{-1}\left( \x - \xfm \right) - \Hl^\top \R^{-1}\left( \Hn(\x) - \y \right).
    \end{split}
\end{equation}
\end{example}

However, there are no guarantees about constraint preservation when considering \cref{eq:vfpito} coupled with \cref{eq:vfpoptimaldrift} and the simplifying assumptions. 
Stochastic diffusion without any special treatment almost always destroys constraints.
Deterministic dynamics with $\boldsymbol{\sigma} = \*0$ do not preserve constraints either.
Thus, drawing inspiration from the ``ODE with invariants'' section from Ascher and Petzold~\cite{Ascher_1998_DAE}, we proceed as follows.
The first idea is to consider \cref{eq:vfpito} with the constraint $\g(\x) = \*0$ as a stochastic differential equation with an invariant.
Next, through the equivalence between an ODE with an invariant and a Hessenberg index-2 DAE (see equations 9.38, 9.39 in Ascher~\cite{Ascher_1998_DAE}), we can reformulate \cref{eq:vfpito} as a Hessenberg index-2 stochastic DAE.
%

%%%%%%%%%%%%%%%%%%%%%%%%%%%%%%%%%%%%%%%%%%%%%%%%%%%%%%%%%%%%%%%%%%%%%%%%
\subsection{Ensemble Variational Fokker Planck Filters with Stabilization (VFPSTAB)}
\label{subsec:VFPSTAB}
%%%%%%%%%%%%%%%%%%%%%%%%%%%%%%%%%%%%%%%%%%%%%%%%%%%%%%%%%%%%%%%%%%%%%%%%

In the VFPSTAB approach, the SDE in \cref{eq:vfpito} is stabilized for the constraint $\g(\x)$.
Following Asher and Petzold~\cite[equation 9.40]{Ascher_1998_DAE}, we rewrite \cref{eq:vfpito} as 
\ifreport
\begin{equation}
    \label{eq:eq1stab}
    \dif \xt = \left( \*F(\syt, \xt, \Pt) - \gamma \gx^\dagger(\xt) \g(\xt) \right) \dif \syt + \boldsymbol{\sigma}(\syt, \xt, \Pt) \dif \*W,
\end{equation}
\else
\begin{equation}
\label{eq:eq1stab}
\begin{split}
    \dif \xt &= \left( \*F(\syt, \xt, \Pt) - \gamma \gx^\dagger(\syt, \x_\syt) \g(\xt) \right) \dif \syt\\
    &+ \boldsymbol{\sigma}(\syt, \xt, \Pt) \dif \*W,
\end{split}
\end{equation}
\fi
where $\gamma > 0$ is a tunable stabilization parameter and $\gx^\dagger = \gx^\top (\gx \gx^\top)^{-1}$ is the Moore-Penrose inverse of the Jacobian $\gx$ \cref{eqn:constraint-Jacobian}.
Note that rather than $\gx^\dagger$, one may use operators such as $\widehat{\gx} (\gx \widehat{\gx})^{-1}$ where $\widehat{\gx}$ is any operator such that $\gx \widehat{\gx}$ is always boundedly invertible~\cite{Ascher_1998_DAE}. 
The choice for $\widehat{\gx}$ defines the projection direction and choosing $\widehat{\gx} = \gx^\top$, as in \cref{eq:eq1stab} results in orthogonal projection.   
VFPSTAB is evolved in time using either Euler-Maruyama~\cite{Kloeden_2011_sdebook} or Rosenbrock-Euler-Maruyama~\cite{Subrahmanya_2023_Ensemble}.

%%%%%%%%%%%%%%%%%%%%%%%%%%%%%%%%%%%%%%%%%%%%%%%%%%%%%%%%%%%%%%%%%%%%%%%%
\subsection{Stochastic Differential Algebraic Equation (VFPDAE)}
\label{subsec:vfpdae}
%%%%%%%%%%%%%%%%%%%%%%%%%%%%%%%%%%%%%%%%%%%%%%%%%%%%%%%%%%%%%%%%%%%%%%%%

In the VFPDAE approach, \cref{eq:vfpito} is equipped with constraint $\g(\x)$ as a Hessenberg index-2 stochastic differential algebraic equation (SDAE).
Formally, by introducing an algebraic variable $\z$ into \cref{eq:vfpito}, we can write
\ifreport
\begin{equation}\label{eq:eq1dae}
\begin{split}
    \dif \xt &= \left(\*F(\syt, \xt, \Pt) - \gx^\top(\xt) \; \z \right)\dif \syt + \boldsymbol{\sigma}(\syt, \xt, \Pt) \dif \*W, \\
    \*0 &= \g(\xt),
\end{split}
\end{equation}
\else
\begin{equation}\label{eq:eq1dae}
\begin{split}
    \dif \xt &= \left(\*F(\syt, \xt, \Pt) - \gx^\top(\xt) \; \z \right)\dif \syt\\
    &+ \boldsymbol{\sigma}(\syt, \xt, \Pt) \dif \*W, \\
    \*0 &= \g(\x_\syt),
\end{split}
\end{equation}
\fi
where $\gx(\xt) = \nabla_{\x} \g(\xt)$. 
The algebraic variable $\z \in \Rspace^{\ncon}$ projects the solution along $\gx^\top$, compensates not just for the deterministic dynamics but also the stochastic Wiener process, which almost surely pushes the states off the constraint manifold.
\begin{remark}\label{rem:remelim}
The constraint $\*0 = \g(\xt + \dif \xt)$, implies the tangent constraint $\*0 = \gx(\xt) \dif \xt$.
Plugging $\dif \xt$ from \cref{eq:eq1dae} into $\*0 = \gx(\xt) \dif \xt$, we get
\ifreport
\begin{equation}
    \*0 = \gx(\xt) \left( \left(\*F(\syt, \xt, \Pt) - \gx^\top(\xt) \; \z \right)\dif \syt + \boldsymbol{\sigma}(\syt, \xt, \Pt) \dif \*W \right),
\end{equation}
\else
\begin{equation}
\begin{split}
    \*0 &= \gx(\xt) \left(\*F(\syt, \xt, \Pt) - \gx^\top(\xt) \; \z \right)\dif \syt\\
    &+ \gx(\xt) \boldsymbol{\sigma}(\syt, \xt, \Pt) \dif \*W,
\end{split}
\end{equation}
\fi
which gives 
\ifreport
\begin{equation}
  \z \,\dif \syt = (\gx \gx^\top)^{-1} \gx (\xt) \left( \*F(\syt, \xt, \Pt) \dif \syt + \boldsymbol{\sigma}(\syt, \xt, \Pt) \dif \*W \right).
\end{equation}        
\else
\begin{equation}
\begin{split}
    \z \, \dif \syt &= (\gx \gx^\top)^{-1} \gx (\xt) \*F(\syt, \xt, \Pt) \dif \syt\\
    &+ \boldsymbol{\sigma}(\syt, \xt, \Pt) \dif \*W .
\end{split}
\end{equation}        
\fi
Plugging $\z \dif \syt$ back in \cref{eq:eq1dae}, we get the SDAE system
\ifreport
\begin{equation}
\begin{split}
    \dif \xt &= (\mathbf{I} -  \gx^\top(\gx \gx^\top)^{-1}\gx(\x_\syt))\cdot (\*F(\syt, \xt, \Pt) \,\dif \syt + \boldsymbol{\sigma}(\syt, \xt, \Pt) \dif \*W ) \\
    \*0 &= \g(\x_\syt).    
\end{split}
\end{equation}
\else
\begin{equation}
\begin{split}
    \dif \xt &= (\mathbf{I} -  \gx^\top(\gx \gx^\top)^{-1}\gx(\x_\syt))\cdot \*F(\syt, \xt, \Pt) \,\dif \syt\\ &+ \boldsymbol{\sigma}(\syt, \xt, \Pt) \dif \*W  \\
    \*0 &= \g(\x_\syt).    
\end{split}
\end{equation}
\fi
The $(\mathbf{I} -  \gx^\top(\gx \gx^\top)^{-1}\gx(\x_\syt))$ term is projecting the stochastic dynamics orthogonally onto the constraint tangent manifold, while eliminating the algebraic variable $\z$.
Note that $\gx \gx^\top$ must always be non-singular for the system to have index-2 uniformly. 
\end{remark}
\begin{remark}
    Compared to VFPDAE (as in \cref{eq:eq1dae}), VFPSTAB (as in \cref{eq:eq1stab}) does not satisfy the constraints exactly, but merely makes the states remain \textit{close} to the constraint manifold.
    However, VFPSTAB is computationally cheaper than VFPDAE due to not having to solve a root-finding problem at each evolution step.
    Both methods improve the preservation of constraints when compared to the standard VFP.
\end{remark}
%

%%%%%%%%%%%%%%%%%%%%%%%%%%%%%%%%%%%%%%%%%%%%%%%%%%%%%%%%%%%%%%%%%%%%%%%%
\subsubsection{Time integration of stochastic differential-algebraic equations}
\label{subsubsec:timeintegration}
%%%%%%%%%%%%%%%%%%%%%%%%%%%%%%%%%%%%%%%%%%%%%%%%%%%%%%%%%%%%%%%%%%%%%%%%

While theory exists for the numerical solution of SDE~\cite{Kloeden_2011_sdebook} and index-2 DAEs~\cite{Ascher_1998_DAE,Wanner_1996_Solving,Hairer_2006_DAE}, very minimal literature exists for index-2 SDAE.
Due to this, we combine ideas from SDE and DAE literature to evolve an index-2 SDAE.
A direct half-explicit solution to the SDAE in \cref{eq:eq1dae} is obtained by taking an Euler-Maruyama step, followed by a projection onto the constraint manifold:
\begin{equation}\label{eq:eq1daedisc1}
\begin{split}
    \tilde{\x}_{1} &= \x_0 + h \*F(\syt_0, \x_0) + \sqrt{h} \boldsymbol{\sigma}(\syt_0, \x_0) \boldsymbol{\xi},\\
    \text{Solve for } \*z_0 &: \g(\tilde{\x}_1 - \gx^\top(\x_0)\, \z_0) = \*0, \\
    \x_1 &= \tilde{\x}_1 - \gx^\top(\x_0)\, \z_0,
\end{split}
\end{equation}
where $h = \syt_1 - \syt_0$.
This method projects the solution along the adjoint constraint manifold evaluated at $\x_0$ (as the constraint does not change over pseudo-time $\syt$). 
Alternatively, one may discretize \cref{eq:eq1dae}
\begin{equation}\label{eq:eq1daedisc2}
\begin{split}
    \tilde{\x}_{1} &= \x_0 + h \*F(\syt_0, \x_0) + \sqrt{h} \boldsymbol{\sigma}(\syt_0, \x_0) \boldsymbol{\xi},\\
    \text{Solve for } \*z_1 &: \g(\tilde{\x}_1 - \gx^\top(\tilde{\x}_1)\, \z_1) = \*0, \\
    \x_1 &= \tilde{\x}_1 - \gx^\top(\tilde{\x}_1)\, \z_1.
\end{split}
\end{equation}
The difference between \cref{eq:eq1daedisc1} and \cref{eq:eq1daedisc2} lies in the evaluation of $\gx$ at $\x_0$ versus $\tilde{\x}_1$. 
This only changes the projection direction $\gx^\top$ to be computed at $\x_1$.
\Cref{eq:eq1daedisc2} can be seen as an evolve and project approach where the solution is evolved by the posterior flow dynamic using Euler-Maruyama~\cite{Kloeden_2011_sdebook} discretization which is then projected onto the constraint manifold.
The flow can also be discretized with the Rosenbrock Euler-Maruyama (discussed in our work \cite{Subrahmanya_2023_Ensemble}), followed by the projection:
\ifreport
\begin{equation}\label{eq:eq1daedisc3}
\begin{split}
    \tilde{\x}_{1} &= \x_0 + h \left( \*I - h\*F_{\x}(\syt_0, \x_0)  \right)^{-1} \*F(\syt_0, \x_0) + \sqrt{h} \boldsymbol{\sigma}(\syt_0, \x) \boldsymbol{\xi},\\
    \text{Solve for } \*z_1 &: \g(\tilde{\x}_1 - \gx^\top(\tilde{\x}_1)\, \z_1) = \*0, \\
    \x_1 &= \tilde{\x}_1 - \gx^\top(\tilde{\x}_1)\, \z_1.
\end{split}
\end{equation}
\else
\begin{equation}\label{eq:eq1daedisc3}
\begin{split}
    \tilde{\x}_{1} &= \x_0 + h \left( \*I - h\*F_{\x}(\syt_0, \x_0)  \right)^{-1} \*F(\syt_0, \x_0)\\
    &+ \sqrt{h} \boldsymbol{\sigma}(\syt_0, \x) \boldsymbol{\xi},\\
    \text{Solve for } \*z_1 &: \g(\tilde{\x}_1 - \gx^\top(\tilde{\x}_1)\, \z_1) = \*0, \\
    \x_1 &= \tilde{\x}_1 - \gx^\top(\tilde{\x}_1)\, \z_1.
\end{split}
\end{equation}
\fi

Note that the rootfinding for $\z_1$ in \cref{eq:eq1daedisc1,eq:eq1daedisc2,eq:eq1daedisc3}  is typically low-dimensional (of dimension $\Rspace^{\ncon}$) compared to the state dimension, i.e. $\ncon \ll \nstate$.
One may also solve for $\tilde{\x}_1$ using any technique from \cite{Crouse_2019_consideration, Crouse_2021_Stiff} and project the states onto the constraint manifold.
For the listed methods \cref{eq:eq1daedisc1,eq:eq1daedisc2,eq:eq1daedisc3}, the solution $\tilde{\x}_1$ will be half-th order accurate ($\!O(\sqrt{h})$).
Since the algebraic variable $\z$ is just the projection vector, we ignore whatever order it may be. 
In our experiments, we found \cref{eq:eq1daedisc2,eq:eq1daedisc3} to give faster results, and these methods have been used for all experiments.

\begin{remark}
    As in \cref{rem:remelim}, we can evolve the projected SDE as
    \ifreport
    \begin{equation}\label{eq:rem3}
        \tilde{\x}_{1} = \x_0 + \left(\*{I} -  \gx^\top(\gx \gx^\top)^{-1}\gx(\x_0)  \right)^{-1} \left(h \*F(\syt_0, \x_0) + \sqrt{h} \boldsymbol{\sigma}(\syt_0, \x) \boldsymbol{\xi}\right),
    \end{equation}
    \else
    \begin{equation}\label{eq:rem3}
    \begin{split}
        \tilde{\x}_{1} &= \x_0 + \left(\*{I} -  \gx^\top(\gx \gx^\top)^{-1}\gx(\x_0)  \right)^{-1}\\
        &\cdot \left(h \*F(\syt_0, \x_0) + \sqrt{h} \boldsymbol{\sigma}(\syt_0, \x) \boldsymbol{\xi}\right),
    \end{split}
    \end{equation}
    \fi
    and project $\tilde{\x}_1$ onto the constraint manifold. 
    In practice, \cref{eq:rem3} took longer wall-time to converge to the analysis, without any qualitative difference from the solutions of \cref{eq:eq1daedisc2,eq:eq1daedisc3}. 
    For this reason, the results from \cref{eq:rem3} are not reported in the paper.
\end{remark}

%%%%%%%%%%%%%%%%%%%%%%%%%%%%%%%%%%%%%%%%%%%%%%%%%%%%%%%%%%%%%%%%%%%%%%%%
\section{Numerical experiments}
\label{sec:expts}
%%%%%%%%%%%%%%%%%%%%%%%%%%%%%%%%%%%%%%%%%%%%%%%%%%%%%%%%%%%%%%%%%%%%%%%%

The VFPSTAB and VFPDAE methodologies are tested on three different test problems, namely:
\begin{enumerate*}[label={(\roman*)}]
    \item the double pendulum,
    \item the Korteweg-de Vries equation, and
    \item the incompressible Navier Stokes via barotropic vorticity equation.
\end{enumerate*}
These results are compared to the results from the standard VFP, the standard ETKF, ETKFP (ETKF with projection as described in~\cref{subsec:proj}), and ETKFA (ETKF augmented with pseudo-observations as described in~\cref{subsec:pseudo}).
For the incompressible Navier-Stokes, we also look at LETKFP and LETKFA which are R-localized~\cite{Reich_2015_book} versions of ETKFP and ETKFA respectively. 

Two different metrics---the spatio-temporal root mean square error (RMSE) of the analysis states and the spatio-temporal root mean square error of the analysis constraints (CRMSE)---are used to evaluate the results.
The RMSE at time $t_k$ is defined as 
\begin{equation}
    \operatorname{RMSE}(k) = \sqrt{\frac{1}{C_1}\sum_{i = \rho}^k \sum_{e = 1}^{\nens} \left\lVert \xai^{[e]} - \xtri \right\rVert_2^2 }, 
\end{equation}
where $C_1 = (k - \rho + 1)\nens \nstate$ where $\nstate$ is the number of states, $\nens$ is the number of ensemble members, and $\rho$ is the spinup steps before which the error statistics are ignored. 
We chose to not use the traditional spatio-temporal RMSE formula $\sqrt{\frac{1}{(k - \rho + 1) \nstate}\sum_{i = \rho}^k \left\lVert \xami - \xtri \right\rVert_2^2 }$, with $\xami = \frac{1}{\nens}\sum_{e = 1}^{\nens}\xai^{[e]}$, since the mean of an ensemble $\xami$ will not lie on the constraint manifold.
In an ideal situation, one would have to compute the mean on the constraint manifold and compute the errors using distances on the constraint manifold.
This approach is not pursued due to the complexities of computing distances on a manifold. 
The CRMSE at time $t_k$ is defined as 
\begin{equation}
    \operatorname{CRMSE}(k) = \sqrt{\frac{1}{C_2} \sum_{i = \rho}^k \sum_{e = 1}^{\nens} \left \lVert \*{E} \cdot  \*{g}(\xai^{[e]}) \right\rVert^2_{2} }, 
\end{equation}
where $C_2 = (k - \rho + 1) \nens \ncon$ and $\*{E}$ is a chosen scaling matrix. 
%

%%%%%%%%%%%%%%%%%%%%%%%%%%%%%%%%%%%%%%%%%%%%%%%%%%%%%%%%%%%%%%%%%%%%%%%%
\subsection{Covariance Regularization}
\label{subsubsec:covreg}
%%%%%%%%%%%%%%%%%%%%%%%%%%%%%%%%%%%%%%%%%%%%%%%%%%%%%%%%%%%%%%%%%%%%%%%%

Covariance estimation from a small ensemble is an important and challenging issue in data assimilation~\cite{Matsuo_2018_Covariance}.  
Various justifiable heuristics such as inflation~\cite{Anderson_1999_MC-implementation}, localization~\cite{Houtekamer_2001}, covariance shrinkage~\cite{Chen_2010_shrinkage}, and modified Cholesky~\cite{Elias_MC_2018}, have been developed to improve the estimates.

For constraint-preserving systems the covariance in full state space is singular, and the empirically estimated ensemble variances are small along directions normal to the constraint manifold, leading to filter divergence. 
In the VFP, the computed drift in \cref{eq:gaussiandrift} requires a linear system solution of the covariance matrix. 
If the covariance is either singular or if its entries are too small, the computed drift will be inaccurate leading to filter divergence over multiple assimilation cycles. 
For this reason, we propose covariance regularization approaches to prevent VFP filter divergence for the dynamics of \cref{eq:eq1stab,eq:eq1dae}.

The estimated covariance is regularized via static covariance shrinkage for the double pendulum.
The estimated covariance $\*P$ is shrunk to $\gamma_{sh} \*P + (1 - \gamma_{sh}) \*I $, with $\gamma_{sh} = 0.01$.
While there are robust methods to compute $\gamma_{sh}$ such as the Ledoit-Wolf and the Rao-Blackwell Ledoit-Wolf methods~\cite{Chen_2010_shrinkage}, we prefer an experimentally tuned static value for quicker computation.

For high-dimensional PDE systems, it is well-known that the square root of the covariance matrix can approximated by the inverse of a spatial Laplacian-like operator~\cite{Stuart_2010_Laplacian,Tan_2013_IP} by its structure (specifically, the inverse of the Laplacian-like operator is diagonally dominant, symmetric, and positive definite). 
The precision matrix estimate for the Kortweg-de Vries and the incompressible Navier-Stokes equations is given by scaling the discrete spatial Laplacian-like of the system with a scalar which is the inverse of the maximum (across the $\nstate$ states) variance of the ensemble. 
For the incompressible Navier-Stokes problem, we assume that the velocities in the $x$ and $y$ dimensions are uncorrelated. 
The specific Laplacian-like approximations are described in \cref{subsec:exptkdv} and \cref{subsec:exptins} respectively.

%%%%%%%%%%%%%%%%%%%%%%%%%%%%%%%%%%%%%%%%%%%%%%%%%%%%%%%%%%%%%%%%%%%%%%%%
\subsection{Double Pendulum}
\label{subsec:exptdp}
%%%%%%%%%%%%%%%%%%%%%%%%%%%%%%%%%%%%%%%%%%%%%%%%%%%%%%%%%%%%%%%%%%%%%%%%

The double pendulum is a constrained mechanical system ~\cite{Hairer_2006_DAE,Ascher_1998_DAE}.
This problem was chosen due to its chaotic dynamical nature and ease of constraint visualization while posing a challenging assimilation problem.
We consider the index-2 DAE formulation of the pendulum dynamics
\begin{equation}
    \frac{\dif^2}{\dif t^2}
    \begin{bmatrix}
        x_1 \\ y_1 \\ x_2 \\ y_2        
    \end{bmatrix} = 
    \begin{bmatrix}
        - \lambda_1 x_1 + \lambda_2 (x_2 - x_1)\\
        - \lambda_1 y_1 + \lambda_2 (y_2 - y_1) - g_c \\
        - \lambda_2 (x_2 - x_1)\\
        - \lambda_2 (y_2 - y_1) - g_c\\  
    \end{bmatrix},
\end{equation}
where $g_c = 9.8$ is the acceleration due to gravity, $x_1$, $y_1$, $x_2$, $y_2$ are the positions of the two pendulums in Cartesian coordinates with rod tensions $\lambda_1, \lambda_2$ being the solution of the linear system
\begin{equation}
    \begin{bmatrix}
        x_1^2 + y_1^2 & -(x_1\delta_{x} + y_1\delta_{y}) \\
        -(x_1\delta_{x} + y_1\delta_{y}) & 2(\delta_{x}^2 + \delta_{y}^2)
    \end{bmatrix} 
    \begin{bmatrix} \lambda_1 \\ \lambda_2 \end{bmatrix} =  
    \begin{bmatrix}
        u_1^2 + v_1^2 - g_c y_1 \\ \delta_{u}^2 + \delta_{v}^2  
    \end{bmatrix},
\end{equation}
where $u_1 = \frac{\dif x_1}{\dif t}$, $v_1 = \frac{\dif y_1}{\dif t}$, $u_2 = \frac{\dif x_2}{\dif t}$, $v_2 = \frac{\dif y_2}{\dif t}$ are the velocities of the two pendulums in the $x$ and $y$ directions, and $\delta_{x} = (x_2 - x_1)$, $\delta_{y} = (y_2 - y_1)$, $\delta_{u} = (u_2 - u_1)$, $\delta_{v} = (v_2 - v_1)$.
The system has $\ncon = 5$ constraints on the lengths of the two rods, velocities of the two rods, and total mechanical energy respectively defined as
\begin{align}
    \begin{bmatrix}
        0 \\ 0 \\ 0 \\ 0 \\ 0
    \end{bmatrix}
    =
    \begin{bmatrix}
        0.5(x_1^2 + y_1^2 - 1) \\
        0.5(\delta_x^2 + \delta_y^2 - 1) \\
        x_1 u_1 + y_1 v_1 \\
        \delta_x \delta_u + \delta_y \delta_v \\
        0.5(u_1^2 + v_1^2 + u_2^2 + v_2^2) + g_c (y_1 + y_2 + 3) - E_0
    \end{bmatrix},
\end{align}
where $E_0$ is the total mechanical energy at the initial condition.

The masses and lengths are assumed to be unity, i.e. $m_1 = m_2 = 1$ and $l_1 = l_2 = 1$.
This index-2 system evolved as a first-order index-2 DAE and using a two-stage, partitioned half explicit Runge-Kutta (PHERK) method~\cite{Murua_1996_Partitioned,Arnold_1998_pherk} for $\Delta t_{ts} = 0.01$.
The coefficients of the second-order PHERK method~\cite{Murua_1996_Partitioned,Arnold_1998_pherk} are 
\begin{equation}
        A = \begin{bmatrix}
            0 & 0 \\ 1 & 0
        \end{bmatrix},
        \tilde{A} = \begin{bmatrix}
            0 & 0 \\ \frac{1}{2} & \frac{1}{2}
        \end{bmatrix},
        c = \tilde{c} = \begin{bmatrix}
            0 \\ 1
        \end{bmatrix},
        b = \begin{bmatrix}
            \frac{1}{2} \\ \frac{1}{2}
        \end{bmatrix}.
\end{equation}
The different matrices $A$ and $\tilde{A}$ are used to solve for the differential and algebraic variables following the work of Murua~\cite{Murua_1996_Partitioned}.

The observation operator $\Hn(\x) = \x$ observes all states with an unbiased, independent Gaussian error $\obserr \sim N(\*0, \*R= 0.1\,\*I_8)$ every $\Delta t = 0.1$ time units.
The problem is assimilated for 5501 steps (each step at $\Delta t = 0.1$) and ignores the statistics from the first $\sigma = 501$ spinup steps. 
A reference system state is chosen to be some position of the double pendulum with zero velocities having a total mechanical energy of $E_0 = 56.1741$.
This reference is then evolved $\nens = 30$ times at a $\Delta t_{samp} = 0.008$ to create a trajectory of states.
The truth and the ensemble are state samples without repetition from this trajectory.

The assimilation methods and their parameters are as follows. 
\begin{enumerate}
    \item ETKF, ETKFA, and ETKFP all use a covariance inflation of $\alpha = 1.08$.
    \item ETKFA uses the augmented error covariance $\R^{\g} = 0.001\*I_5$.
    This value was hand-tuned to balance between the RMSE and CRMSE. 
    \item All VFP variants---VFP, VFPSTAB, VFPDAE---have the diffusion parameter 
    \begin{equation*}
        \boldsymbol{\sigma} = 0.001 \operatorname{diag}\left(\begin{bmatrix}
        2 & 2 & 20 & 20 & 2 & 2 & 20 & 20
    \end{bmatrix}\right).
    \end{equation*}
    This was tuned to be a scaling of the climatological (or auto) covariance of the double pendulum system.
    The VFP system was evolved to the posterior in pseudo-time $\syt$ using Euler-Maruyama with a $\Delta \syt = 0.001$.
    The evolution was terminated when the absolute change in the ensemble mean across pseudo-time steps was under a specified tolerance, here $1e^{-6}$.
    The covariances $\CovPt$, and $\CovPf$ were estimated using covariance shrinkage as discussed in \Cref{subsubsec:covreg}. 
    \item VFPSTAB uses a hand-tuned stabilization constant of $\gamma = 30$  (to balance between RMSE and CRMSE).
\end{enumerate}
\begin{figure*}[!ht]
    \centering
    \begin{subfigure}[b]{0.49\textwidth}
        \centering
        \includegraphics[width=\linewidth]{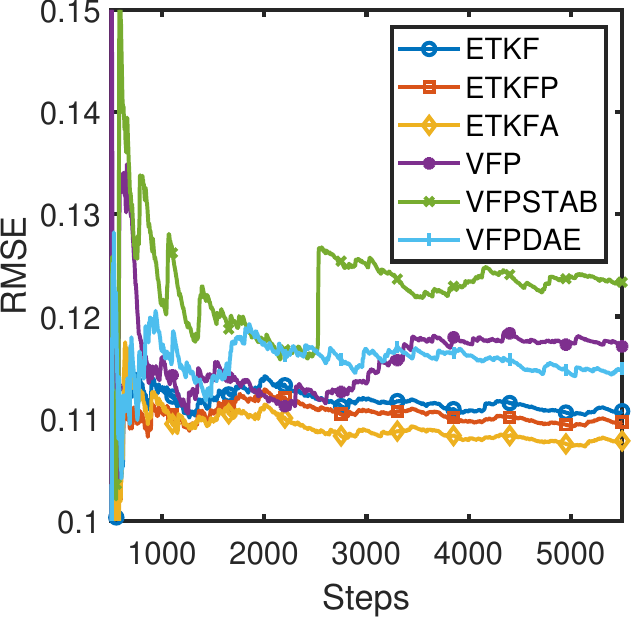}
        \caption{RMSE for various methods. 
        }
        \label{fig:expt_dp1}
    \end{subfigure}
    \hfill
    \begin{subfigure}[b]{0.49\textwidth}
        \centering
        \includegraphics[width=\linewidth]{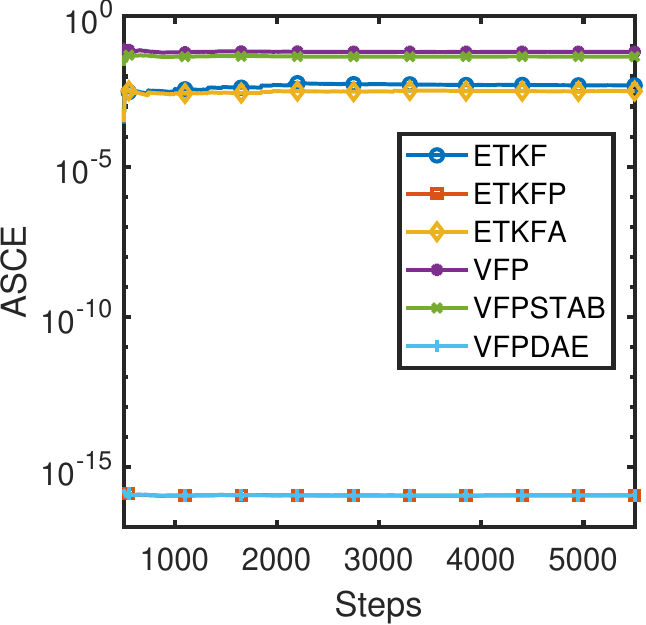}
        \caption{CRMSE for various methods.}
        \label{fig:expt_dp2}
    \end{subfigure}
    \caption{Results for the double pendulum experiments. }
    \label{fig:expt_dp}
\end{figure*}

\Cref{fig:expt_dp1} compares the spatio-temporal root mean squared errors of various flavors of ETKF and VFP. 
For this problem, all the ETKF variants behave slightly better than all the VFP variants, with VFPDAE coming very close to plain ETKF.
For the ETKF methods, constraint preservation (ETKFA and ETKFP) shows a slightly lower RMSE when compared to plain ETKF.
\Cref{fig:expt_dp2} shows the CRMSE for each method with the scaling $\*E = \operatorname{diag}\left(\begin{bmatrix}
    1 & 1 & 1 & 1 & \frac{1}{E_0}
\end{bmatrix}\right)$.
The asymptotic CRMSE values of ETKF, ETKFA, VFP, and VFPSTAB are 0.005, 0.003, 0.065, and 0.046 respectively. 
As expected, ETKFP and VFPDAE preserve the constraints to machine precision.
ETKFA and ETKF have the next lowest CRMSE.
Due to constraint augmentation, ETKFA has a slightly lower CRMSE than ETKF.
Similarly, the CRMSE of VFPSTAB and VFP look identical, but VFPSTAB has a slightly lower CRMSE due to stabilization.
%

%%%%%%%%%%%%%%%%%%%%%%%%%%%%%%%%%%%%%%%%%%%%%%%%%%%%%%%%%%%%%%%%%%%%%%%%
\subsection{Korteweg de Vries}
\label{subsec:exptkdv}
%%%%%%%%%%%%%%%%%%%%%%%%%%%%%%%%%%%%%%%%%%%%%%%%%%%%%%%%%%%%%%%%%%%%%%%%

The Korteweg-de Vries PDE represents non-linear, non-dissipative waves on shallow water surfaces~\cite{DeJager_2006_KDVorigin}. 
While this PDE has an infinite set of integrable constraints~\cite{Miura_2003_KdV} on the states, our experiments consider the first three only.
The equation is written as: 
\begin{equation}
    \frac{\partial \x}{\partial t} + 3 \frac{\partial (\x^2)}{\partial x}  + \frac{\partial^3 \x}{\partial x^3} = 0,
\end{equation}
where $\x$ is a non-dimensional wave-height displacement, and $t$ is the non-dimensional time.
The mass, momentum, and energy integral invariants are respectively given by
\begin{align}
    \phi_1(\x) &= \int_{\Omega_x} \x \mathrm{d} x, \nonumber\\
    \phi_2(\x) &= \int_{\Omega_x} \frac{1}{2} \x^2 \mathrm{d} x, \\
    \phi_3(\x) &= \int_{\Omega_x} \left( \frac{\x^3}{3} - \left( \frac{\partial \x}{\partial x} \right)^2 \right) \mathrm{d} x. \nonumber
\end{align}
Thus, the $\ncon = 3$ constraints, are 
\begin{equation}\label{eq:kdvcon}
    \*0 = \g(t, \x) = \phi(\x) - \phi(\x_0).
\end{equation}

The PDE is spatially discretized using the second-order central finite difference in the periodic domain $x \in [-10, 10]$ with $\nstate = 100$ points.
This semi-discretized PDE is evolved using the implicit midpoint scheme (with $\Delta t_{ts} = 0.01$) which is known to preserve the invariants~\cite{Ascher_2005_Symplectic}.
Following \cite{Dutykh_2013_kdv}, the third constraint in \cref{eq:kdvcon} is discretized with the first-order forward finite difference. 

At time intervals of $\Delta t = 0.01$, every fourth state is observed as $\Hn(\x) = \x(4:4:100) $ with an unbiased Gaussian observation error with covariance $\*R = 0.2\,\*I_{25}$.
The observation covariance corresponds to 10\% of the climatological (or auto) covariance.
Assimilation is done for 2201 steps, and the reported results ignore the statistics from the first $\rho = 401$ spinup steps.
The initial true state $\xtro = 6 \operatorname{sech}^2(x)$ corresponds to the two soliton problem.
The constraints for this are $\phi(\x_0) = \begin{bmatrix} 12 & 48 & -211.3815\end{bmatrix}^\top$.
The initial ensemble with $\nens = 10$ samples is created by perturbing $\xtro$ with diagonal random noise and projecting this onto the constraint manifold.

We consider the following hand-tuned data assimilation methods: 
\begin{enumerate}
    \item ETKF, ETKFA, and ETKFP all use an inflation of $\alpha = 1.04$ at every assimilation step. 
    \item ETKFA uses an augmented error covariance $\R^{\g} = 10\,\*{I}_5$ for the constraints.
    \item VFPSTAB, and VFPDAE use a hand-tuned diffusion of $\boldsymbol{\sigma} = 1e^{-4}(\Delta + 1e^{-3}\*{I}_{\nstate})^{-1}$ where $\Delta$ is the discrete Laplacian.
    As stated in \cref{subsubsec:covreg}, the covariances $\Pt^{-1}$, and $\Pf^{-1}$ are estimated by multiplying $(\Delta + 1e^{-3}\*{I}_{\nstate})^2$ with the maximum variance of the intermediate and forecast ensembles respectively.
    The term $1e^{-3}\*{I}_{\nstate}$ is added to the discrete Laplacian as the discrete Laplacian is low-rank due to the periodic domain. 
    Additionally, perturbed observations~\cite{Houtekamer_1998a} were used in the VFP methods to allow for a larger analysis spread with perturbations being drawn from $\!N( \*0, 0.05^2 \R)$.
    This is because increasing the diffusion results in irregular solutions that do not lie on the constraint manifold.
    \item 
    The VFPSTAB and VFPDAE system are solved in pseudo-time $\syt$ using Rosenbrock-Euler-Maruyama with a $\Delta \syt = 0.01$.
    The evolution is terminated when the absolute change in the ensemble mean across pseudo-time steps falls under a specified tolerance, here $1e^{-4}$.
\end{enumerate}
The LETKF~\cite{Hunt_2007_4DLETKF} (and LETKFP, LETKFA) results are not reported as they are similar to the unlocalized ETKF (and ETKFP, ETKFA) results. 

\begin{figure*}[!ht]
    \centering
    \begin{subfigure}[b]{0.49\textwidth}
        \centering
        \includegraphics[width=\linewidth]{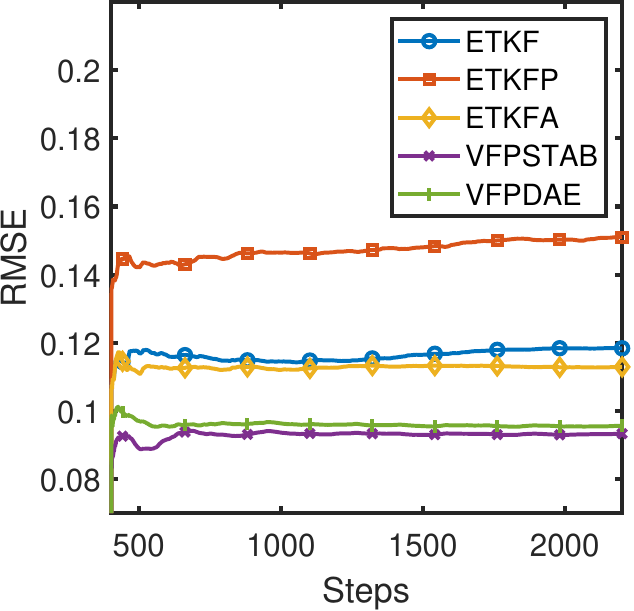}
        \caption{RMSE for various methods.}
        \label{fig:expt_kdv1}
    \end{subfigure}
    \hfill
    \begin{subfigure}[b]{0.49\textwidth}
        \centering
        \includegraphics[width=\linewidth]{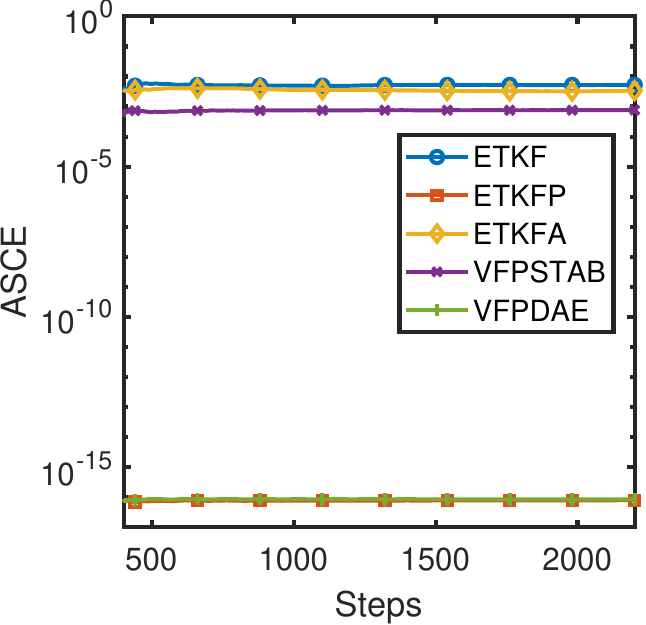}
        \caption{CRMSE for various methods.}
        \label{fig:expt_kdv2}
    \end{subfigure}
    \caption{Results for the Korteweg-deVries experiments where the VFP methods outperform the ETKF methods.}
    \label{fig:expt_kdv}
\end{figure*}

\Cref{fig:expt_kdv1} shows the spatio-temporal RMSE in time for various analysis schemes.
The VFP methods (VFPDAE and VFPSTAB) perform better than the ETKF variants.
Between the ETKF methods, ETKFA seems to perform better than ETKF and ETKFP for the RMSE. 
This might mean that the ETKFP is attaining a worse solution on the manifold, while ETKFA is obtaining a relatively good solution on the manifold compared to vanilla ETKF.
VFPDAE outperforms all these methods.
\Cref{fig:expt_kdv2} shows the CRMSE values for the different assimilation methods with the scaling matrix $\*E = \*I_3$.
The asymptotic CRMSE values of ETKF, ETKFA, and VFPSTAB are 0.005, 0.003, and 0.0007 respectively. 
The asymptotic CRMSE values of ETKFP and VFPDAE are at machine epsilon.
ETKF is worse at preserving constraints than ETKFA, both of which are worse than VFPSTAB.
VFPSTAB is worse than ETKFP and VFPDAE which preserve the constraints exactly.
%

%%%%%%%%%%%%%%%%%%%%%%%%%%%%%%%%%%%%%%%%%%%%%%%%%%%%%%%%%%%%%%%%%%%%%%%%
\subsection{Incompressible Navier-Stokes}
\label{subsec:exptins}
%%%%%%%%%%%%%%%%%%%%%%%%%%%%%%%%%%%%%%%%%%%%%%%%%%%%%%%%%%%%%%%%%%%%%%%%

The incompressible Navier-Stokes equations (in the vorticity-streamfunction form) are used to model large-scale, oceanic flows in shallow basins~\cite{San_2011_QG}. 
Without external forcing, this system must conserve the mass, while also conserving the total energy and enstrophy~\cite{Arakawa_1972_Design,Kacimi_2013_Arakawa}. 
This problem was chosen since it has both local and global constraints.

The equations are 
\begin{equation}
\label{eq:QG}
    \begin{split}
        \frac{\partial \omega}{\partial t} &= - J(\psi, \omega) + Re^{-1} \Delta \omega + Ro^{-1} \left(\frac{\partial \psi}{\partial x} + F \right),\\
        \psi &= - \Delta^{-1} \omega, 
    \end{split}
\end{equation}
with the vorticity $\omega$, streamfunction $\psi$, velocities $(u, v) = \left( \frac{\partial \psi}{\partial y}, -\frac{\partial \psi}{\partial x} \right)$, the Laplacian $\Delta = \frac{\partial^2}{ \partial x^2} + \frac{\partial^2}{ \partial y^2}$ , the Reynolds number $Re = 450$, the Rossby number $Ro = 0.0036$, the Jacobian $J(\psi, \omega) \equiv \frac{\partial \psi}{\partial y} \frac{\partial \omega}{\partial x} - \frac{\partial \psi}{\partial x} \frac{\partial \omega}{\partial y}$, and symmetric double gyre forcing $F = \sin{(\pi y)}$.
The streamfunction and vorticity boundaries on $x$ and $y$ are homogeneous Dirichlet.
The settings are chosen from San et. al~\cite{San_2011_QG} to produce four gyre circulation over time.

The experiment setup involves solving the vorticity equation in time and then computing the primitive variables (velocities $(u, v)$) as the forecast. 
After constraint preserving assimilation on the velocities, the velocities are converted back to the vorticity form as $\omega = \frac{\partial v}{\partial x} - \frac{\partial u}{\partial y}$ with appropriate boundary conditions and marched forward in time.
This approach allows the incompressible Navier-Stokes to be solved as an ODE without computational difficulties.

Firstly, the primitive variables must be incompressible due to the conservation of mass, i.e. $\frac{\partial u}{\partial x} + \frac{\partial v}{\partial y} = 0$.
Next, additional constraints are constructed to preserve the energy and enstrophy for each ensemble member before and after filtering.
This means that each ensemble member will have its own, independent constraint.
The incompressibility, energy, and enstrophy constraints are respectively
\begin{align}
    \frac{\partial u_a}{\partial x} + \frac{\partial v_a}{\partial y} &= 0 , \nonumber\\
    \frac{1}{2} \int_\Omega \left( u_a^2 + v_a^2 \right) \dif x \; \dif y &= \frac{1}{2} \int_\Omega \left( u_f^2 + v_f^2 \right) \dif x \; \dif y ,\\
    \frac{1}{2} \int_\Omega \left( \frac{\partial v_a}{\partial x} - \frac{\partial u_a}{\partial y} \right)^2 \dif x\;  \dif y &= \frac{1}{2} \int_\Omega \left( \frac{\partial v_f}{\partial x} - \frac{\partial u_f}{\partial y} \right)^2 \dif x\;  \dif y. \nonumber
\end{align}
The problem is discretized on the domain $x \in [0, 1]$ and $y \in [-1, 1]$ with central finite differences using $64$ and $129$ equidistant points in each dimension.
This makes $\nstate = 16512$ from having 2 variables ($u$ and $v$) on the grid. 
This discretization results in $8256$ local constraints due to incompressibility, and $2$ global constraints for the energy and enstrophy, making $\ncon = 8258$.
The Arakawa Jacobian~\cite{Kacimi_2013_Arakawa} with the RK3~\cite{Gottlieb_1998_RK3} method is used to evolve the semidiscretized PDE with a timestep of $\Delta_{ts} = 0.0109/40$.

The system is observed every $\Delta t = 0.0109$ units equivalent to one real day with this non-dimensionalization~\cite{San_2011_QG}. 
The velocities $u$ and $v$ are observed at $16 \times 16 = 256$ equally spaced grid points in the interior of the domain, with a Gaussian observation error covariance of $\R = 400 \*I_{512}$. 
This observation error covariance is approximately 1\% of the climatological (or auto) covariance of the velocities. 
The assimilation was done for 881 steps neglecting the statistics from the first 81 steps as spinup.
The true initial condition is generated by evolving a smooth random vorticity field for 100 days. 
Random sampling around the truth generates the initial ensemble with $\nens = 10$.
The following filtering methods are considered for this problem: 
\begin{enumerate}
    \item ETKF and ETKFP with hand-tuned inflations of $\alpha = 1.1$ and $\alpha = 1.05$, respectively.
    \item LETKF, and LETKFP, with an inflation of $\alpha = 1.05$. 
    Here, LETKF and LETKFP use the Gaspari-Cohn localization~\cite{Gaspari_1999_correlation,Houtekamer_2001} with a hand-tuned radius of $r_{loc} = 0.559$ units.
    The distance used to compute the localization coefficients is the spatial  Euclidean distance.
    \item VFPDAE with diffusion $\boldsymbol{\sigma} =  \operatorname{blkdiag}\left(\begin{bmatrix}
        \Delta & \Delta
    \end{bmatrix}\right)^{-1}$, where $\Delta$ is the two-dimensional discrete Laplacian on the grid.
    Note that this assumes that the correlation between $u$ and $v$ on the same grid point is $0$.
    From the authors' experience, it was observed that the velocities in the $x$ and $y$ directions were uncorrelated, making this heuristic meaningful.
    As discussed in \cref{subsubsec:covreg}, $\Pt^{-1}$, and $\Pf^{-1}$ are computed by scaling $\operatorname{blkdiag}\left(\begin{bmatrix}
        \Delta^2 & \Delta^2
    \end{bmatrix}\right)$ with the maximum variance of the intermediate and forecast ensembles respectively.
    Also, to accelerate convergence, we initialize the VFPDAE solution with a projected solution from ETKF.
    The VFPDAE system initialized at the projected ETKF solution was evolved to the analysis with the flow dynamic for 10 steps using Euler-Maruyama with a $\Delta \syt = 1e^{-6}$.
    Due to the cost of this method, the pseudo-time evolution was stopped after 10 evolution steps.
    \item Augmented observations could not be used for the ETKF and LETKF as the constraint is different for each ensemble member. 
\end{enumerate}
\begin{figure}
    \centering
    \ifreport
    \includegraphics[width=0.5\linewidth]{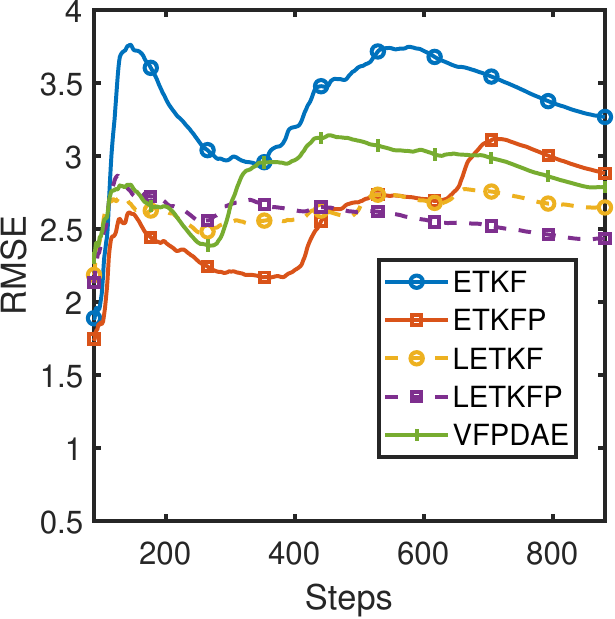}
    \else
    \includegraphics[width=\linewidth]{figs/BVE/bvermses.pdf}
    \fi
    \caption{Incompressible Navier-Stokes: RMSE for different methods.}
    \label{fig:expt_bve1}
\end{figure}
The results in \cref{fig:expt_bve1} show that ETKF and ETKFP methods perform worse than their localized variants LETKF and LETKFP.
Our method VFPDAE performs better than ETKF and ETKFP, but worse than LETKF and LETKFP.
Even with 10 small steps of size $1e^{-5}$, VFPDAE performs better than its initialization method ETKFP.
Perhaps more evolution pseudo-time could give better results but this was not done due to the high computational cost.

\begin{figure*}[!ht]
    \centering
    \begin{subfigure}[b]{0.49\textwidth}
        \centering
        \includegraphics[width=\linewidth]{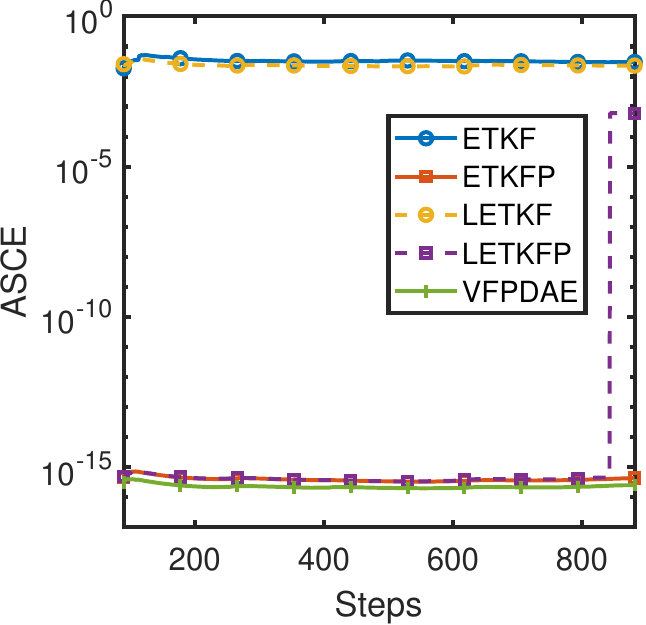}
        \caption{Enstrophy CRMSE for different methods.}
        \label{fig:expt_bve2}
    \end{subfigure}
    \hfill
    \begin{subfigure}[b]{0.49\textwidth}
        \centering
        \includegraphics[width=\linewidth]{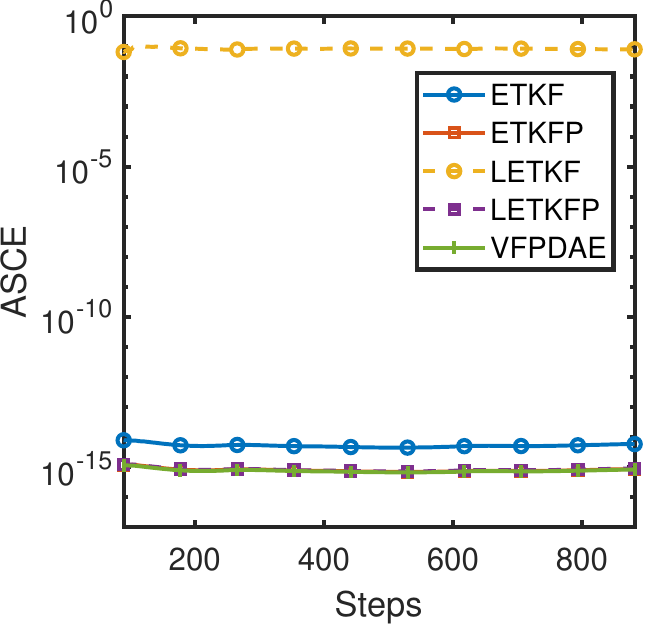}
        \caption{Divergence CRMSE for different methods.}
        \label{fig:expt_bve3}
    \end{subfigure}
    \caption{CRMSE results for the Navier-Stokes experiments. }
    \label{fig:expt_bvea}
\end{figure*}
\Cref{fig:expt_bve2} chooses an $\*E = \operatorname{diag} \left(\begin{bmatrix}
    0 & \dots & 0 & 1
\end{bmatrix}\right))$ for looking only at the enstrophy. 
The asymptotic CRMSE values for the ETKF, LETKF, and LETKFP, are 0.03, 0.02, and 6e-4 respectively.  
VFPDAE and ETKFP preserve the enstrophy to the numerical limit, while ETKF and LETKF do not.
LETKFP preserves enstrophy to the numerical limit at all times except step 840, which makes its CRMSE go up.
This happens because the projection method fails to project the LETKF solution onto the constraint manifold. 

A more interesting situation is looking at the spatially averaged CRMSE for only the divergence constraint across the states with $\*E = \operatorname{diag}\left(\begin{bmatrix}
    1 & \dots & 1 & 0 & 0
\end{bmatrix}\right)$ in \Cref{fig:expt_bve3} .
LETKF violates the divergence constraint by a huge margin, with the final time CRMSE being 0.08. 
It is important to notice that ETKF preserves it naturally (around an asymptotic CRMSE of ten times machine epsilon) due to it being a linear transform of the ensemble as also observed by \cite{Janjic_2014_Preserve}.
ETKFP, LETKFP, and VFPDAE preserve the divergence constraint to the numerical limit.

%%%%%%%%%%%%%%%%%%%%%%%%%%%%%%%%%%%%%%%%%%%%%%%%%%%%%%%%%%%%%%%%%%%%%%%%
\section{Conclusions}
\label{sec:conc}
%%%%%%%%%%%%%%%%%%%%%%%%%%%%%%%%%%%%%%%%%%%%%%%%%%%%%%%%%%%%%%%%%%%%%%%%

The VFP framework is extended to incorporate general non-linear equality constraints.
The two formulated methods---VFPDAE and VFPSTAB---respect the constraint manifold throughout the evolution to the posterior.
To preserve constraints \textit{exactly}, one must use the VFPDAE method which requires a somewhat expensive implicit-explicit time-stepping method to evolve the dynamics. 
VFPSTAB allows one to make computational cost gains at the cost of \textit{inexactly} preserving the constraints.
Both these methods require covariance regularization as the constraint manifold can restrict the particles to have small variance, leading to filter divergence over time. 
While not tested, non-linear inequality and box constraints can be incorporated by drawing on ideas from constrained optimization theory such as using an ``active-set'' of constraints in the flow, enabling VFP to handle more generic constraints.
The same constraint-preserving ideas can be extended to smoothing over a trajectory of observations.
Future work will include work on both smoothing and inequality constraints. 

\backmatter

%\bmhead{Supplementary information}
\bmhead{Acknowledgements}

This work was supported by DOE through award ASCR DE-SC0021313, by NSF through award  CDS\&E--MSS 1953113, and by the Computational Science Laboratory at Virginia Tech. The authors also thank Dr. Traian Iliescu from the Mathematics Department at Virginia Tech for helpful discussions.

\section*{Declarations}
\bmhead{Competing interests} All authors certify that they have no affiliations with or involvement in any organization or entity with any financial interest or non-financial interest in the subject matter or materials discussed in this manuscript.

\bmhead{Data availability} Code will be made available on request by email.

%\begin{appendices}
%\end{appendices}

\bibliography{Bib/references,Bib/data_assim_particle,Bib/data_assim_general,Bib/data_assim_models,Bib/data_assim_kalman,Bib/sandu}


\begin{thebibliography}{61}
\ifx \bisbn   \undefined \def \bisbn  #1{ISBN #1}\fi
\ifx \binits  \undefined \def \binits#1{#1}\fi
\ifx \bauthor  \undefined \def \bauthor#1{#1}\fi
\ifx \batitle  \undefined \def \batitle#1{#1}\fi
\ifx \bjtitle  \undefined \def \bjtitle#1{#1}\fi
\ifx \bvolume  \undefined \def \bvolume#1{\textbf{#1}}\fi
\ifx \byear  \undefined \def \byear#1{#1}\fi
\ifx \bissue  \undefined \def \bissue#1{#1}\fi
\ifx \bfpage  \undefined \def \bfpage#1{#1}\fi
\ifx \blpage  \undefined \def \blpage #1{#1}\fi
\ifx \burl  \undefined \def \burl#1{\textsf{#1}}\fi
\ifx \doiurl  \undefined \def \doiurl#1{\url{https://doi.org/#1}}\fi
\ifx \betal  \undefined \def \betal{\textit{et al.}}\fi
\ifx \binstitute  \undefined \def \binstitute#1{#1}\fi
\ifx \binstitutionaled  \undefined \def \binstitutionaled#1{#1}\fi
\ifx \bctitle  \undefined \def \bctitle#1{#1}\fi
\ifx \beditor  \undefined \def \beditor#1{#1}\fi
\ifx \bpublisher  \undefined \def \bpublisher#1{#1}\fi
\ifx \bbtitle  \undefined \def \bbtitle#1{#1}\fi
\ifx \bedition  \undefined \def \bedition#1{#1}\fi
\ifx \bseriesno  \undefined \def \bseriesno#1{#1}\fi
\ifx \blocation  \undefined \def \blocation#1{#1}\fi
\ifx \bsertitle  \undefined \def \bsertitle#1{#1}\fi
\ifx \bsnm \undefined \def \bsnm#1{#1}\fi
\ifx \bsuffix \undefined \def \bsuffix#1{#1}\fi
\ifx \bparticle \undefined \def \bparticle#1{#1}\fi
\ifx \barticle \undefined \def \barticle#1{#1}\fi
\bibcommenthead
\ifx \bconfdate \undefined \def \bconfdate #1{#1}\fi
\ifx \botherref \undefined \def \botherref #1{#1}\fi
\ifx \url \undefined \def \url#1{\textsf{#1}}\fi
\ifx \bchapter \undefined \def \bchapter#1{#1}\fi
\ifx \bbook \undefined \def \bbook#1{#1}\fi
\ifx \bcomment \undefined \def \bcomment#1{#1}\fi
\ifx \oauthor \undefined \def \oauthor#1{#1}\fi
\ifx \citeauthoryear \undefined \def \citeauthoryear#1{#1}\fi
\ifx \endbibitem  \undefined \def \endbibitem {}\fi
\ifx \bconflocation  \undefined \def \bconflocation#1{#1}\fi
\ifx \arxivurl  \undefined \def \arxivurl#1{\textsf{#1}}\fi
\csname PreBibitemsHook\endcsname

%%% 1
\bibitem[\protect\citeauthoryear{Asch et~al.}{2016}]{Asch_2016_book}
\begin{bbook}
\bauthor{\bsnm{Asch}, \binits{M.}},
\bauthor{\bsnm{Bocquet}, \binits{M.}},
\bauthor{\bsnm{Nodet}, \binits{M.}}:
\bbtitle{Data Assimilation}.
\bpublisher{Society for Industrial and Applied Mathematics},
\blocation{Philadelphia, PA}
(\byear{2016}).
\doiurl{10.1137/1.9781611974546} .
\burl{https://epubs.siam.org/doi/abs/10.1137/1.9781611974546}
\end{bbook}
\endbibitem

%%% 2
\bibitem[\protect\citeauthoryear{Reich and Cotter}{2015}]{Reich_2015_book}
\begin{bbook}
\bauthor{\bsnm{Reich}, \binits{S.}},
\bauthor{\bsnm{Cotter}, \binits{C.}}:
\bbtitle{Probabilistic Forecasting and {Bayesian} Data Assimilation}.
\bpublisher{Cambridge University Press},
\blocation{Cambridge, UK}
(\byear{2015})
\end{bbook}
\endbibitem

%%% 3
\bibitem[\protect\citeauthoryear{Evensen et~al.}{2022}]{Evensen_2022_book}
\begin{bbook}
\bauthor{\bsnm{Evensen}, \binits{G.}},
\bauthor{\bsnm{Vossepoel}, \binits{F.C.}},
\bauthor{\bsnm{Leeuwen}, \binits{P.J.}}:
\bbtitle{Data Assimilation Fundamentals: A Unified Formulation of the State and
  Parameter Estimation Problem}.
\bpublisher{Springer},
\blocation{Cham}
(\byear{2022}).
\doiurl{10.1007/978-3-030-96709-3} .
\burl{https://doi.org/10.1007/978-3-030-96709-3}
\end{bbook}
\endbibitem

%%% 4
\bibitem[\protect\citeauthoryear{Kalman}{1960}]{Kalman_1960}
\begin{barticle}
\bauthor{\bsnm{Kalman}, \binits{R.E.}}:
\batitle{A new approach to linear filtering and prediction problems}.
\bjtitle{Transaction of the ASME- Journal of Basic Engineering}
\bvolume{82},
\bfpage{35}--\blpage{45}
(\byear{1960})
\end{barticle}
\endbibitem

%%% 5
\bibitem[\protect\citeauthoryear{Evensen}{1994}]{Evensen_1994}
\begin{barticle}
\bauthor{\bsnm{Evensen}, \binits{G.}}:
\batitle{Sequential data assimilation with a nonlinear quasi-geostrophic model
  using {Monte Carlo} methods to forcast error statistics}.
\bjtitle{Journal of Geophysical Research}
\bvolume{99}(\bissue{C5}),
\bfpage{10143}--\blpage{10162}
(\byear{1994})
\end{barticle}
\endbibitem

%%% 6
\bibitem[\protect\citeauthoryear{Bishop et~al.}{2001}]{Bishop_2001_ETKF}
\begin{barticle}
\bauthor{\bsnm{Bishop}, \binits{C.}},
\bauthor{\bsnm{Etherton}, \binits{B.J.}},
\bauthor{\bsnm{Majumdar}, \binits{S.}}:
\batitle{{Adaptive sampling with the Ensemble Transform Kalman Filter. Part I:
  Theoretical Aspects}}.
\bjtitle{Monthly Weather Review}
\bvolume{129},
\bfpage{420}--\blpage{436}
(\byear{2001})
\end{barticle}
\endbibitem

%%% 7
\bibitem[\protect\citeauthoryear{Anderson}{2001}]{Anderson_2001_EAKF}
\begin{barticle}
\bauthor{\bsnm{Anderson}, \binits{J.L.}}:
\batitle{An {e}nsemble {a}djustment {K}alman {f}ilter for data assimilation}.
\bjtitle{Monthly Weather Review}
\bvolume{129},
\bfpage{2884}--\blpage{2903}
(\byear{2001})
\end{barticle}
\endbibitem

%%% 8
\bibitem[\protect\citeauthoryear{Hunt et~al.}{2007}]{Hunt_2007_4DLETKF}
\begin{barticle}
\bauthor{\bsnm{Hunt}, \binits{B.R.}},
\bauthor{\bsnm{Kostelich}, \binits{E.J.}},
\bauthor{\bsnm{Szunyogh}, \binits{I.}}:
\batitle{Efficient data assimilation for spatiotemporal chaos: A local ensemble
  transform {K}alman filter}.
\bjtitle{Physica D: Nonlinear Phenomena}
\bvolume{230}(\bissue{1}),
\bfpage{112}--\blpage{126}
(\byear{2007})
\doiurl{10.1016/j.physd.2006.11.008} .
\bcomment{Data Assimilation}
\end{barticle}
\endbibitem

%%% 9
\bibitem[\protect\citeauthoryear{van Leeuwen}{2009}]{vanLeeuwen_2009_PF-review}
\begin{barticle}
\bauthor{\bsnm{Leeuwen}, \binits{P.J.}}:
\batitle{Particle filtering in geophysical systems}.
\bjtitle{Monthly Weather Review}
\bvolume{137},
\bfpage{4089}--\blpage{4114}
(\byear{2009})
\end{barticle}
\endbibitem

%%% 10
\bibitem[\protect\citeauthoryear{van Leeuwen
  et~al.}{2019}]{vanLeeuwen_2019_PF-review}
\begin{barticle}
\bauthor{\bsnm{Leeuwen}, \binits{P.J.}},
\bauthor{\bsnm{K{\"u}nsch}, \binits{H.R.}},
\bauthor{\bsnm{Nerger}, \binits{L.}},
\bauthor{\bsnm{Potthast}, \binits{R.}},
\bauthor{\bsnm{Reich}, \binits{S.}}:
\batitle{Particle filters for high-dimensional geoscience applications: A
  review}.
\bjtitle{Quarterly Journal of the Royal Meteorological Society}
\bvolume{145}(\bissue{723}),
\bfpage{2335}--\blpage{2365}
(\byear{2019})
\doiurl{10.1002/qj.3551}
{\href{https://arxiv.org/abs/https://rmets.onlinelibrary.wiley.com/doi/pdf/10.1002/qj.3551}{{https://rmets.onlinelibrary.wiley.com/doi/pdf/10.1002/qj.3551}}}
\end{barticle}
\endbibitem

%%% 11
\bibitem[\protect\citeauthoryear{Pulido and {van
  Leeuwen}}{2019}]{Pulido_2019_mapping-PF}
\begin{barticle}
\bauthor{\bsnm{Pulido}, \binits{M.}},
\bauthor{\bsnm{{van Leeuwen}}, \binits{P.J.}}:
\batitle{Sequential {Monte Carlo} with kernel embedded mappings: The mapping
  particle filter}.
\bjtitle{Journal of Computational Physics}
\bvolume{396},
\bfpage{400}--\blpage{415}
(\byear{2019})
\doiurl{10.1016/j.jcp.2019.06.060}
\end{barticle}
\endbibitem

%%% 12
\bibitem[\protect\citeauthoryear{Hu and van Leeuwen}{2021}]{Hu_2020_mapping-PF}
\begin{barticle}
\bauthor{\bsnm{Hu}, \binits{C.-C.}},
\bauthor{\bsnm{Leeuwen}, \binits{P.J.}}:
\batitle{A particle flow filter for high-dimensional system applications}.
\bjtitle{Quarterly Journal of the Royal Meteorological Society}
\bvolume{147}(\bissue{737}),
\bfpage{2352}--\blpage{2374}
(\byear{2021})
\doiurl{10.1002/qj.4028}
{\href{https://arxiv.org/abs/https://rmets.onlinelibrary.wiley.com/doi/pdf/10.1002/qj.4028}{{https://rmets.onlinelibrary.wiley.com/doi/pdf/10.1002/qj.4028}}}
\end{barticle}
\endbibitem

%%% 13
\bibitem[\protect\citeauthoryear{Pathiraja and
  Reich}{2019}]{Reich_2019_discrete-gradients}
\begin{barticle}
\bauthor{\bsnm{Pathiraja}, \binits{S.}},
\bauthor{\bsnm{Reich}, \binits{S.}}:
\batitle{Discrete gradients for computational {B}ayesian inference}.
\bjtitle{Journal of Computational Dynamics}
\bvolume{6}(\bissue{2158-2491\_2019\_2\_385}),
\bfpage{385}
(\byear{2019})
\doiurl{10.3934/jcd.2019019}
\end{barticle}
\endbibitem

%%% 14
\bibitem[\protect\citeauthoryear{Reich and
  Weissmann}{2021}]{Reich_2021_FokkerPlanck}
\begin{barticle}
\bauthor{\bsnm{Reich}, \binits{S.}},
\bauthor{\bsnm{Weissmann}, \binits{S.}}:
\batitle{{Fokker--Planck Particle Systems for Bayesian Inference: Computational
  Approaches}}.
\bjtitle{SIAM/ASA Journal on Uncertainty Quantification}
\bvolume{9}(\bissue{2}),
\bfpage{446}--\blpage{482}
(\byear{2021})
\doiurl{10.1137/19M1303162}
{\href{https://arxiv.org/abs/https://doi.org/10.1137/19M1303162}{{https://doi.org/10.1137/19M1303162}}}
\end{barticle}
\endbibitem

%%% 15
\bibitem[\protect\citeauthoryear{Garbuno-Inigo
  et~al.}{2020}]{Stuart_2020_gradient-EnKF}
\begin{barticle}
\bauthor{\bsnm{Garbuno-Inigo}, \binits{A.}},
\bauthor{\bsnm{Hoffmann}, \binits{F.}},
\bauthor{\bsnm{Li}, \binits{W.}},
\bauthor{\bsnm{Stuart}, \binits{A.M.}}:
\batitle{{Interacting Langevin Diffusions: Gradient Structure and Ensemble
  Kalman Sampler}}.
\bjtitle{SIAM Journal on Applied Dynamical Systems}
\bvolume{19}(\bissue{1}),
\bfpage{412}--\blpage{441}
(\byear{2020})
\doiurl{10.1137/19M1251655}
{\href{https://arxiv.org/abs/https://doi.org/10.1137/19M1251655}{{https://doi.org/10.1137/19M1251655}}}
\end{barticle}
\endbibitem

%%% 16
\bibitem[\protect\citeauthoryear{Crouse and
  Lewis}{2019}]{Crouse_2019_consideration}
\begin{botherref}
\oauthor{\bsnm{Crouse}, \binits{D.F.}},
\oauthor{\bsnm{Lewis}, \binits{C.}}:
Consideration of particle flow filter implementations and biases.
Naval Research Laboratory Memo,
1--17
(2019)
\end{botherref}
\endbibitem

%%% 17
\bibitem[\protect\citeauthoryear{Subrahmanya
  et~al.}{2023}]{Subrahmanya_2023_Ensemble}
\begin{botherref}
\oauthor{\bsnm{Subrahmanya}, \binits{A.N.}},
\oauthor{\bsnm{Popov}, \binits{A.A.}},
\oauthor{\bsnm{Sandu}, \binits{A.}}:
Ensemble Variational Fokker-Planck Methods for Data Assimilation
(2023)
\end{botherref}
\endbibitem

%%% 18
\bibitem[\protect\citeauthoryear{Hastie et~al.}{2009}]{Hastie_2001_statsbook}
\begin{bbook}
\bauthor{\bsnm{Hastie}, \binits{T.}},
\bauthor{\bsnm{Tibshirani}, \binits{R.}},
\bauthor{\bsnm{Friedman}, \binits{J.}}:
\bbtitle{The Elements of Statistical Learning: Data Mining, Inference, and
  Prediction}.
\bpublisher{Springer},
\blocation{New York, NY}
(\byear{2009}).
\doiurl{10.1007/978-0-387-84858-7} .
\burl{https://doi.org/10.1007/978-0-387-84858-7}
\end{bbook}
\endbibitem

%%% 19
\bibitem[\protect\citeauthoryear{Arakawa}{1972}]{Arakawa_1972_Design}
\begin{botherref}
\oauthor{\bsnm{Arakawa}, \binits{A.}}:
Design of the ucla general circulation model.
Technical report,
University of California, Los Angeles
(1972)
\end{botherref}
\endbibitem

%%% 20
\bibitem[\protect\citeauthoryear{Zeng and Janjić}{2016}]{Janjic_2016_LETKF}
\begin{barticle}
\bauthor{\bsnm{Zeng}, \binits{Y.}},
\bauthor{\bsnm{Janjić}, \binits{T.}}:
\batitle{Study of conservation laws with the local ensemble transform {K}alman
  filter}.
\bjtitle{Quarterly Journal of the Royal Meteorological Society}
\bvolume{142}(\bissue{699}),
\bfpage{2359}--\blpage{2372}
(\byear{2016})
\doiurl{10.1002/qj.2829}
{\href{https://arxiv.org/abs/https://rmets.onlinelibrary.wiley.com/doi/pdf/10.1002/qj.2829}{{https://rmets.onlinelibrary.wiley.com/doi/pdf/10.1002/qj.2829}}}
\end{barticle}
\endbibitem

%%% 21
\bibitem[\protect\citeauthoryear{Janjić et~al.}{2014}]{Janjic_2014_Preserve}
\begin{barticle}
\bauthor{\bsnm{Janjić}, \binits{T.}},
\bauthor{\bsnm{McLaughlin}, \binits{D.}},
\bauthor{\bsnm{Cohn}, \binits{S.E.}},
\bauthor{\bsnm{Verlaan}, \binits{M.}}:
\batitle{Conservation of mass and preservation of positivity with ensemble-type
  {K}alman filter algorithms}.
\bjtitle{Monthly Weather Review}
\bvolume{142}(\bissue{2}),
\bfpage{755}--\blpage{773}
(\byear{2014})
\doiurl{10.1175/MWR-D-13-00056.1}
\end{barticle}
\endbibitem

%%% 22
\bibitem[\protect\citeauthoryear{Massicotte
  et~al.}{1995}]{Massicotte_1995_Positivity}
\begin{barticle}
\bauthor{\bsnm{Massicotte}, \binits{D.}},
\bauthor{\bsnm{Morawski}, \binits{R.Z.}},
\bauthor{\bsnm{Barwicz}, \binits{A.}}:
\batitle{Incorporation of a positivity constraint into a {K}alman-filter-based
  algorithm for correction of spectrometric data}.
\bjtitle{IEEE Transactions on Instrumentation and Measurement}
\bvolume{44}(\bissue{1}),
\bfpage{2}--\blpage{7}
(\byear{1995})
\doiurl{10.1109/19.368111}
\end{barticle}
\endbibitem

%%% 23
\bibitem[\protect\citeauthoryear{Simon and Chia}{2002}]{Simon_2002_Equality}
\begin{barticle}
\bauthor{\bsnm{Simon}, \binits{D.}},
\bauthor{\bsnm{Chia}, \binits{T.L.}}:
\batitle{{K}alman filtering with state equality constraints}.
\bjtitle{IEEE Transactions on Aerospace and Electronic Systems}
\bvolume{38}(\bissue{1}),
\bfpage{128}--\blpage{136}
(\byear{2002})
\doiurl{10.1109/7.993234}
\end{barticle}
\endbibitem

%%% 24
\bibitem[\protect\citeauthoryear{Julier and LaViola}{2007}]{Julier_2007_Nonlin}
\begin{barticle}
\bauthor{\bsnm{Julier}, \binits{S.J.}},
\bauthor{\bsnm{LaViola}, \binits{J.J.}}:
\batitle{On {K}alman filtering with nonlinear equality constraints}.
\bjtitle{IEEE Transactions on Signal Processing}
\bvolume{55}(\bissue{6}),
\bfpage{2774}--\blpage{2784}
(\byear{2007})
\doiurl{10.1109/TSP.2007.893949}
\end{barticle}
\endbibitem

%%% 25
\bibitem[\protect\citeauthoryear{Gupta and
  Hauser}{2007}]{Gupta_2007_Constraint}
\begin{botherref}
\oauthor{\bsnm{Gupta}, \binits{N.}},
\oauthor{\bsnm{Hauser}, \binits{R.}}:
{K}alman Filtering with Equality and Inequality State Constraints
(2007)
\end{botherref}
\endbibitem

%%% 26
\bibitem[\protect\citeauthoryear{Sircoulomb
  et~al.}{2008}]{Sircoulomb_2008_Ineq}
\begin{bchapter}
\bauthor{\bsnm{Sircoulomb}, \binits{V.}},
\bauthor{\bsnm{Hoblos}, \binits{G.}},
\bauthor{\bsnm{Chafouk}, \binits{H.}},
\bauthor{\bsnm{Ragot}, \binits{J.}}:
\bctitle{State estimation under nonlinear state inequality constraints. a
  tracking application}.
In: \bbtitle{2008 16th Mediterranean Conference on Control and Automation},
pp. \bfpage{1669}--\blpage{1674}
(\byear{2008}).
\doiurl{10.1109/MED.2008.4602024}
\end{bchapter}
\endbibitem

%%% 27
\bibitem[\protect\citeauthoryear{Zanetti et~al.}{2009}]{Zanetti_2009_norm}
\begin{barticle}
\bauthor{\bsnm{Zanetti}, \binits{R.}},
\bauthor{\bsnm{Majji}, \binits{M.}},
\bauthor{\bsnm{Bishop}, \binits{R.H.}},
\bauthor{\bsnm{Mortari}, \binits{D.}}:
\batitle{Norm-constrained {K}alman filtering}.
\bjtitle{Journal of Guidance, Control, and Dynamics}
\bvolume{32}(\bissue{5}),
\bfpage{1458}--\blpage{1465}
(\byear{2009})
\doiurl{10.2514/1.43119}
{\href{https://arxiv.org/abs/https://doi.org/10.2514/1.43119}{{https://doi.org/10.2514/1.43119}}}
\end{barticle}
\endbibitem

%%% 28
\bibitem[\protect\citeauthoryear{Prakash et~al.}{2010}]{Prakash_2010_Ensemble}
\begin{barticle}
\bauthor{\bsnm{Prakash}, \binits{J.}},
\bauthor{\bsnm{Patwardhan}, \binits{S.C.}},
\bauthor{\bsnm{Shah}, \binits{S.L.}}:
\batitle{Constrained nonlinear state estimation using ensemble {K}alman
  filters}.
\bjtitle{Industrial \& Engineering Chemistry Research}
\bvolume{49}(\bissue{5}),
\bfpage{2242}--\blpage{2253}
(\byear{2010})
\doiurl{10.1021/ie900197s}
{\href{https://arxiv.org/abs/https://doi.org/10.1021/ie900197s}{{https://doi.org/10.1021/ie900197s}}}
\end{barticle}
\endbibitem

%%% 29
\bibitem[\protect\citeauthoryear{Bavdekar
  et~al.}{2013}]{Bavdekar_2013_Constrained}
\begin{bchapter}
\bauthor{\bsnm{Bavdekar}, \binits{V.A.}},
\bauthor{\bsnm{Prakash}, \binits{J.}},
\bauthor{\bsnm{Shah}, \binits{S.L.}},
\bauthor{\bsnm{Gopaluni}, \binits{R.B.}}:
\bctitle{Constrained dual ensemble {K}alman filter for state and parameter
  estimation}.
In: \bbtitle{2013 American Control Conference},
pp. \bfpage{3093}--\blpage{3098}
(\byear{2013}).
\doiurl{10.1109/ACC.2013.6580306}
\end{bchapter}
\endbibitem

%%% 30
\bibitem[\protect\citeauthoryear{Kacimi et~al.}{2013}]{Kacimi_2013_Arakawa}
\begin{botherref}
\oauthor{\bsnm{Kacimi}, \binits{A.}},
\oauthor{\bsnm{Aliziane}, \binits{T.}},
\oauthor{\bsnm{Khouider}, \binits{B.}}:
The {Arakawa Jacobian} method and a fourth-order essentially nonoscillatory
  scheme for the beta-plane barotropic equations.
International Journal of Numerical Analysis \& Modeling
\textbf{10}(3)
(2013)
\end{botherref}
\endbibitem

%%% 31
\bibitem[\protect\citeauthoryear{Zeng et~al.}{2017}]{Janjic_2017_Preserve}
\begin{barticle}
\bauthor{\bsnm{Zeng}, \binits{Y.}},
\bauthor{\bsnm{Janjić}, \binits{T.}},
\bauthor{\bsnm{Ruckstuhl}, \binits{Y.}},
\bauthor{\bsnm{Verlaan}, \binits{M.}}:
\batitle{Ensemble-type {K}alman filter algorithm conserving mass, total energy
  and enstrophy}.
\bjtitle{Quarterly Journal of the Royal Meteorological Society}
\bvolume{143}(\bissue{708}),
\bfpage{2902}--\blpage{2914}
(\byear{2017})
\doiurl{10.1002/qj.3142}
{\href{https://arxiv.org/abs/https://rmets.onlinelibrary.wiley.com/doi/pdf/10.1002/qj.3142}{{https://rmets.onlinelibrary.wiley.com/doi/pdf/10.1002/qj.3142}}}
\end{barticle}
\endbibitem

%%% 32
\bibitem[\protect\citeauthoryear{Li et~al.}{2019}]{Li_2019_Constrained}
\begin{barticle}
\bauthor{\bsnm{Li}, \binits{R.}},
\bauthor{\bsnm{Jan}, \binits{N.M.}},
\bauthor{\bsnm{Huang}, \binits{B.}},
\bauthor{\bsnm{Prasad}, \binits{V.}}:
\batitle{Constrained ensemble {K}alman filter based on {K}ullback--{L}eibler
  divergence}.
\bjtitle{Journal of Process Control}
\bvolume{81},
\bfpage{150}--\blpage{161}
(\byear{2019})
\end{barticle}
\endbibitem

%%% 33
\bibitem[\protect\citeauthoryear{Albers et~al.}{2019}]{Albers_2019_Constraint}
\begin{barticle}
\bauthor{\bsnm{Albers}, \binits{D.J.}},
\bauthor{\bsnm{Blancquart}, \binits{P.-A.}},
\bauthor{\bsnm{Levine}, \binits{M.E.}},
\bauthor{\bsnm{Seylabi}, \binits{E.E.}},
\bauthor{\bsnm{Stuart}, \binits{A.}}:
\batitle{Ensemble {K}alman methods with constraints}.
\bjtitle{Inverse Problems}
\bvolume{35}(\bissue{9}),
\bfpage{095007}
(\byear{2019})
\doiurl{10.1088/1361-6420/ab1c09}
\end{barticle}
\endbibitem

%%% 34
\bibitem[\protect\citeauthoryear{Schillings and
  Stuart}{2017}]{Stuart_2017_EnKF-inversion}
\begin{barticle}
\bauthor{\bsnm{Schillings}, \binits{C.}},
\bauthor{\bsnm{Stuart}, \binits{A.M.}}:
\batitle{Analysis of the {Ensemble Kalman Filter for Inverse Problems}}.
\bjtitle{SIAM Journal on Numerical Analysis}
\bvolume{55}(\bissue{3}),
\bfpage{1264}--\blpage{1290}
(\byear{2017})
\doiurl{10.1137/16M105959X}
{\href{https://arxiv.org/abs/https://doi.org/10.1137/16M105959X}{{https://doi.org/10.1137/16M105959X}}}
\end{barticle}
\endbibitem

%%% 35
\bibitem[\protect\citeauthoryear{Chada et~al.}{2019}]{Neil_2019_Box}
\begin{barticle}
\bauthor{\bsnm{Chada}, \binits{N.K.}},
\bauthor{\bsnm{Schillings}, \binits{C.}},
\bauthor{\bsnm{Weissmann}, \binits{S.}}:
\batitle{On the incorporation of box-constraints for ensemble {K}alman
  inversion}.
\bjtitle{Foundations of Data Science}
\bvolume{1}(\bissue{4}),
\bfpage{433}--\blpage{456}
(\byear{2019})
\doiurl{10.3934/fods.2019018}
\end{barticle}
\endbibitem

%%% 36
\bibitem[\protect\citeauthoryear{Bertsekas}{1982}]{Bertsekas_1982_PN}
\begin{barticle}
\bauthor{\bsnm{Bertsekas}, \binits{D.P.}}:
\batitle{Projected {Newton} methods for optimization problems with simple
  constraints}.
\bjtitle{SIAM Journal on control and Optimization}
\bvolume{20}(\bissue{2}),
\bfpage{221}--\blpage{246}
(\byear{1982})
\end{barticle}
\endbibitem

%%% 37
\bibitem[\protect\citeauthoryear{Amor et~al.}{2018}]{Amor_2018_Constrained}
\begin{botherref}
\oauthor{\bsnm{Amor}, \binits{N.}},
\oauthor{\bsnm{Rasool}, \binits{G.}},
\oauthor{\bsnm{Bouaynaya}, \binits{N.C.}}:
Constrained state estimation--a review.
arXiv preprint arXiv:1807.03463
(2018)
\end{botherref}
\endbibitem

%%% 38
\bibitem[\protect\citeauthoryear{Ascher and Petzold}{1998}]{Ascher_1998_DAE}
\begin{bbook}
\bauthor{\bsnm{Ascher}, \binits{U.M.}},
\bauthor{\bsnm{Petzold}, \binits{L.R.}}:
\bbtitle{Computer Methods for Ordinary Differential Equations and
  Differential-Algebraic Equations}
vol. \bseriesno{61}.
\bpublisher{SIAM},
\blocation{Philadelpha, PA}
(\byear{1998}).
\doiurl{10.1137/1.9781611971392} .
\burl{https://epubs.siam.org/doi/book/10.1137/1.9781611971392}
\end{bbook}
\endbibitem

%%% 39
\bibitem[\protect\citeauthoryear{Anco and
  Cheviakov}{2020}]{Anco_2020_Conservation}
\begin{barticle}
\bauthor{\bsnm{Anco}, \binits{S.C.}},
\bauthor{\bsnm{Cheviakov}, \binits{A.F.}}:
\batitle{On the different types of global and local conservation laws for
  partial differential equations in three spatial dimensions: Review and recent
  developments}.
\bjtitle{International Journal of Non-Linear Mechanics}
\bvolume{126},
\bfpage{103569}
(\byear{2020})
\doiurl{10.1016/j.ijnonlinmec.2020.103569}
\end{barticle}
\endbibitem

%%% 40
\bibitem[\protect\citeauthoryear{Garbuno-Inigo
  et~al.}{2020}]{Inigo_2020_Langevin}
\begin{barticle}
\bauthor{\bsnm{Garbuno-Inigo}, \binits{A.}},
\bauthor{\bsnm{N\"{u}sken}, \binits{N.}},
\bauthor{\bsnm{Reich}, \binits{S.}}:
\batitle{{Affine Invariant Interacting Langevin Dynamics for Bayesian
  Inference}}.
\bjtitle{SIAM Journal on Applied Dynamical Systems}
\bvolume{19}(\bissue{3}),
\bfpage{1633}--\blpage{1658}
(\byear{2020})
\doiurl{10.1137/19M1304891}
{\href{https://arxiv.org/abs/https://doi.org/10.1137/19M1304891}{{https://doi.org/10.1137/19M1304891}}}
\end{barticle}
\endbibitem

%%% 41
\bibitem[\protect\citeauthoryear{Kloeden and
  Platen}{1992}]{Kloeden_2011_sdebook}
\begin{bbook}
\bauthor{\bsnm{Kloeden}, \binits{P.E.}},
\bauthor{\bsnm{Platen}, \binits{E.}}:
\bbtitle{Numerical Solution of Stochastic Differential Equations}.
\bsertitle{Stochastic Modelling and Applied Probability}.
\bpublisher{Springer},
\blocation{Berlin, Heidelberg}
(\byear{1992}).
\doiurl{10.1007/978-3-662-12616-5_2} .
\burl{https://doi.org/10.1007/978-3-662-12616-5_2}
\end{bbook}
\endbibitem

%%% 42
\bibitem[\protect\citeauthoryear{Wanner and Hairer}{1996}]{Wanner_1996_Solving}
\begin{bbook}
\bauthor{\bsnm{Wanner}, \binits{G.}},
\bauthor{\bsnm{Hairer}, \binits{E.}}:
\bbtitle{Solving Ordinary Differential Equations II: Stiff and
  Differential-Algebraic Problems}
vol. \bseriesno{375}.
\bpublisher{Springer},
\blocation{Berlin, Heidelberg}
(\byear{1996}).
\doiurl{10.1007/978-3-642-05221-7} .
\burl{https://doi.org/10.1007/978-3-642-05221-7}
\end{bbook}
\endbibitem

%%% 43
\bibitem[\protect\citeauthoryear{Hairer et~al.}{2006}]{Hairer_2006_DAE}
\begin{bbook}
\bauthor{\bsnm{Hairer}, \binits{E.}},
\bauthor{\bsnm{Lubich}, \binits{C.}},
\bauthor{\bsnm{Roche}, \binits{M.}}:
\bbtitle{The Numerical Solution of Differential-Algebraic Systems by
  Runge-Kutta Methods}
vol. \bseriesno{1409}.
\bpublisher{Springer},
\blocation{Berlin, Heidelberg}
(\byear{2006})
\end{bbook}
\endbibitem

%%% 44
\bibitem[\protect\citeauthoryear{Crouse}{2021}]{Crouse_2021_Stiff}
\begin{botherref}
\oauthor{\bsnm{Crouse}, \binits{D.F.}}:
Particle flow solutions avoiding stiff integration.
NRL Technical Report No. NRU5340/FR-2021/1
(2021)
\end{botherref}
\endbibitem

%%% 45
\bibitem[\protect\citeauthoryear{Matsuo}{2018}]{Matsuo_2018_Covariance}
\begin{barticle}
\bauthor{\bsnm{Matsuo}, \binits{T.}}:
\batitle{Covariance modeling in applications of data assimilation to
  high-dimensional earth and geospace systems}.
\bjtitle{Journal of Physics: Conference Series}
\bvolume{1036}(\bissue{1}),
\bfpage{012019}
(\byear{2018})
\doiurl{10.1088/1742-6596/1036/1/012019}
\end{barticle}
\endbibitem

%%% 46
\bibitem[\protect\citeauthoryear{Anderson and
  Anderson}{1999}]{Anderson_1999_MC-implementation}
\begin{barticle}
\bauthor{\bsnm{Anderson}, \binits{J.L.}},
\bauthor{\bsnm{Anderson}, \binits{S.L.}}:
\batitle{A {M}onte {C}arlo implementation of the nonlinear filtering problem to
  produce ensemble assimilations and forecasts}.
\bjtitle{Mon. Wea. Rev.}
\bvolume{127},
\bfpage{2741}--\blpage{2785}
(\byear{1999})
\end{barticle}
\endbibitem

%%% 47
\bibitem[\protect\citeauthoryear{Houtekamer and
  Mitchell}{2001}]{Houtekamer_2001}
\begin{barticle}
\bauthor{\bsnm{Houtekamer}, \binits{P.L.}},
\bauthor{\bsnm{Mitchell}, \binits{H.L.}}:
\batitle{{A sequential ensemble Kalman filter for atmospheric data assimilation
  }}.
\bjtitle{Monthly Weather Review}
\bvolume{129},
\bfpage{123}--\blpage{137}
(\byear{2001})
\end{barticle}
\endbibitem

%%% 48
\bibitem[\protect\citeauthoryear{{Chen} et~al.}{2010}]{Chen_2010_shrinkage}
\begin{barticle}
\bauthor{\bsnm{{Chen}}, \binits{Y.}},
\bauthor{\bsnm{{Wiesel}}, \binits{A.}},
\bauthor{\bsnm{{Eldar}}, \binits{Y.C.}},
\bauthor{\bsnm{{Hero}}, \binits{A.O.}}:
\batitle{{Shrinkage Algorithms for MMSE Covariance Estimation}}.
\bjtitle{IEEE Transactions on Signal Processing}
\bvolume{58}(\bissue{10}),
\bfpage{5016}--\blpage{5029}
(\byear{2010})
\doiurl{10.1109/TSP.2010.2053029}
\end{barticle}
\endbibitem

%%% 49
\bibitem[\protect\citeauthoryear{Nino-Ruiz et~al.}{2018}]{Elias_MC_2018}
\begin{barticle}
\bauthor{\bsnm{Nino-Ruiz}, \binits{E.D.}},
\bauthor{\bsnm{Sandu}, \binits{A.}},
\bauthor{\bsnm{Deng}, \binits{X.}}:
\batitle{An ensemble {K}alman filter implementation based on modified cholesky
  decomposition for inverse covariance matrix estimation}.
\bjtitle{SIAM Journal on Scientific Computing}
\bvolume{40}(\bissue{2}),
\bfpage{867}--\blpage{886}
(\byear{2018})
\doiurl{10.1137/16M1097031}
{\href{https://arxiv.org/abs/https://doi.org/10.1137/16M1097031}{{https://doi.org/10.1137/16M1097031}}}
\end{barticle}
\endbibitem

%%% 50
\bibitem[\protect\citeauthoryear{Stuart}{2010}]{Stuart_2010_Laplacian}
\begin{barticle}
\bauthor{\bsnm{Stuart}, \binits{A.M.}}:
\batitle{Inverse problems: A {Bayesian} perspective}.
\bjtitle{Acta Numerica}
\bvolume{19},
\bfpage{451}--\blpage{559}
(\byear{2010})
\doiurl{10.1017/S0962492910000061}
\end{barticle}
\endbibitem

%%% 51
\bibitem[\protect\citeauthoryear{Bui-Thanh et~al.}{2013}]{Tan_2013_IP}
\begin{barticle}
\bauthor{\bsnm{Bui-Thanh}, \binits{T.}},
\bauthor{\bsnm{Ghattas}, \binits{O.}},
\bauthor{\bsnm{Martin}, \binits{J.}},
\bauthor{\bsnm{Stadler}, \binits{G.}}:
\batitle{A computational framework for infinite-dimensional bayesian inverse
  problems part i: The linearized case, with application to global seismic
  inversion}.
\bjtitle{SIAM Journal on Scientific Computing}
\bvolume{35}(\bissue{6}),
\bfpage{2494}--\blpage{2523}
(\byear{2013})
\doiurl{10.1137/12089586X}
{\href{https://arxiv.org/abs/https://doi.org/10.1137/12089586X}{{https://doi.org/10.1137/12089586X}}}
\end{barticle}
\endbibitem

%%% 52
\bibitem[\protect\citeauthoryear{Murua}{1996}]{Murua_1996_Partitioned}
\begin{botherref}
\oauthor{\bsnm{Murua}, \binits{A.}}:
Partitioned {Runge-Kutta} methods for semi-explicit differential-algebraic
  systems of index 2.
Technical Report EHU-KZAA-IKT-196
(1996)
\end{botherref}
\endbibitem

%%% 53
\bibitem[\protect\citeauthoryear{Arnold}{1998}]{Arnold_1998_pherk}
\begin{barticle}
\bauthor{\bsnm{Arnold}, \binits{M.}}:
\batitle{Half-explicit {Runge-Kutta} methods with explicit stages for
  differential-algebraic systems of index 2}.
\bjtitle{BIT Numerical Mathematics}
\bvolume{38}(\bissue{3}),
\bfpage{415}--\blpage{438}
(\byear{1998})
\doiurl{10.1007/BF02510252}
\end{barticle}
\endbibitem

%%% 54
\bibitem[\protect\citeauthoryear{De~Jager}{2006}]{DeJager_2006_KDVorigin}
\begin{botherref}
\oauthor{\bsnm{De~Jager}, \binits{E.}}:
On the origin of the {Korteweg-de Vries} equation.
arXiv preprint math/0602661
(2006)
\end{botherref}
\endbibitem

%%% 55
\bibitem[\protect\citeauthoryear{Miura et~al.}{2003}]{Miura_2003_KdV}
\begin{barticle}
\bauthor{\bsnm{Miura}, \binits{R.M.}},
\bauthor{\bsnm{Gardner}, \binits{C.S.}},
\bauthor{\bsnm{Kruskal}, \binits{M.D.}}:
\batitle{{Korteweg‐de Vries Equation and Generalizations. II. Existence of
  Conservation Laws and Constants of Motion}}.
\bjtitle{Journal of Mathematical Physics}
\bvolume{9}(\bissue{8}),
\bfpage{1204}--\blpage{1209}
(\byear{2003})
\doiurl{10.1063/1.1664701}
{\href{https://arxiv.org/abs/https://pubs.aip.org/aip/jmp/article-pdf/9/8/1204/8183110/1204\_1\_online.pdf}{{https://pubs.aip.org/aip/jmp/article-pdf/9/8/1204/8183110/1204\_1\_online.pdf}}}
\end{barticle}
\endbibitem

%%% 56
\bibitem[\protect\citeauthoryear{Ascher and
  McLachlan}{2005}]{Ascher_2005_Symplectic}
\begin{barticle}
\bauthor{\bsnm{Ascher}, \binits{U.M.}},
\bauthor{\bsnm{McLachlan}, \binits{R.I.}}:
\batitle{On symplectic and multisymplectic schemes for the {KdV} equation}.
\bjtitle{Journal of Scientific Computing}
\bvolume{25},
\bfpage{83}--\blpage{104}
(\byear{2005})
\end{barticle}
\endbibitem

%%% 57
\bibitem[\protect\citeauthoryear{Dutykh et~al.}{2013}]{Dutykh_2013_kdv}
\begin{barticle}
\bauthor{\bsnm{Dutykh}, \binits{D.}},
\bauthor{\bsnm{Chhay}, \binits{M.}},
\bauthor{\bsnm{Fedele}, \binits{F.}}:
\batitle{Geometric numerical schemes for the {KdV} equation}.
\bjtitle{Computational Mathematics and Mathematical Physics}
\bvolume{53}(\bissue{2}),
\bfpage{221}--\blpage{236}
(\byear{2013})
\doiurl{10.1134/S0965542513020103}
\end{barticle}
\endbibitem

%%% 58
\bibitem[\protect\citeauthoryear{Houtekamer and
  Mitchell}{1998}]{Houtekamer_1998a}
\begin{barticle}
\bauthor{\bsnm{Houtekamer}, \binits{P.L.}},
\bauthor{\bsnm{Mitchell}, \binits{H.L.}}:
\batitle{{Data assimilation using an ensemble Kalman filter technique }}.
\bjtitle{Monthly Weather Review}
\bvolume{126},
\bfpage{796}--\blpage{811}
(\byear{1998})
\end{barticle}
\endbibitem

%%% 59
\bibitem[\protect\citeauthoryear{San et~al.}{2011}]{San_2011_QG}
\begin{barticle}
\bauthor{\bsnm{San}, \binits{O.}},
\bauthor{\bsnm{Staples}, \binits{A.E.}},
\bauthor{\bsnm{Wang}, \binits{Z.}},
\bauthor{\bsnm{Iliescu}, \binits{T.}}:
\batitle{Approximate deconvolution large eddy simulation of a barotropic ocean
  circulation model}.
\bjtitle{Ocean Modelling}
\bvolume{40}(\bissue{2}),
\bfpage{120}--\blpage{132}
(\byear{2011})
\end{barticle}
\endbibitem

%%% 60
\bibitem[\protect\citeauthoryear{Gottlieb and Shu}{1998}]{Gottlieb_1998_RK3}
\begin{barticle}
\bauthor{\bsnm{Gottlieb}, \binits{S.}},
\bauthor{\bsnm{Shu}, \binits{C.-W.}}:
\batitle{Total variation diminishing {Runge-Kutta} schemes}.
\bjtitle{Mathematics of computation}
\bvolume{67}(\bissue{221}),
\bfpage{73}--\blpage{85}
(\byear{1998})
\end{barticle}
\endbibitem

%%% 61
\bibitem[\protect\citeauthoryear{Gaspari and
  Cohn}{1999}]{Gaspari_1999_correlation}
\begin{barticle}
\bauthor{\bsnm{Gaspari}, \binits{G.}},
\bauthor{\bsnm{Cohn}, \binits{S.E.}}:
\batitle{Construction of correlation functions in two and three dimensions}.
\bjtitle{Quarterly Journal of the Royal Meteorological Society}
\bvolume{125},
\bfpage{723}--\blpage{757}
(\byear{1999})
\end{barticle}
\endbibitem

\end{thebibliography}

\end{document}